\documentclass[12pt, a4]{article}

\usepackage{amsmath,amssymb, amsfonts}

\usepackage{amscd}

\title{\bf Axiomatic Rejection for the Propositional Fragment of Le\'{s}niewski's Ontology}
\author{Takao Inou\'{e}, Arata Ishimoto\footnote{Deceased. This paper is dedicated to the memory of Professor Emeritus Arata Ishimoto. This paper is cited in Ishimoto \cite{Ishimoto1996}.} and Mitsunori Kobayashi}
\date{August 15, 2021}

\usepackage{tikz}
\usetikzlibrary{trees}

\newtheorem{theorem}{Theorem}
\newtheorem{lemma}{Lemma}

\newtheorem{corollary}{Corollary} 
\newtheorem{definition}{Definition}

\newtheorem{proposition}{Proposition}

\begin{document}

\maketitle

\begin{abstract}
A Hilbert-type axiomatic rejection $\mathbf{HAR}$ for the propositional fragment $\mathbf{L_1}$ of
Le\'{s}niewski's ontology is proposed. Also a Gentzen-type axiomatic rejection $\mathbf{GAR}$ 
of $\mathbf{L_1}$ is 
proposed. Models for $\mathbf{L_1}$ are introduced. By axiomatic rejection, Ishimoto's embedding theorem 
will be proved. One of our main theorems is:

\noindent \textsc{Theorem}  \rm (Main Theorem) \it 

$\vdash_T A \enspace  \Longleftrightarrow \enspace \vdash_H \enspace A$

$\enspace \enspace  \enspace  \enspace  \enspace \enspace  \Longleftrightarrow 
\enspace TA \enspace \mbox{is valid in first-order predicate logic with equality}$

$\enspace \enspace  \enspace  \enspace  \enspace \enspace  \Longleftrightarrow 
\enspace not \dashv_H A.$ 

\noindent where $\vdash_T A$ means that $A$ is provable in the tableau method of $\mathbf{L_1}$, 
while $\vdash_H A$ means that $A$ is provable in the Hilbert-type $\mathbf{L_1}$. \rm
\smallskip

In the last section, as the chracterization theorem, we shall show: 

\noindent \textsc{Theorem} \rm (Characterization Theorem) \it 
The following statements are equivalent:

$(1)$ The Cut elimination theorem for the Gentzen-type $\mathbf{L_1}$ $($i.e. tableau method$)$ holds,

$(2)$  No Hintikka formula of the form $A_1 \vee A_2 \vee \cdots \vee A_n$ $(n \geq 1)$ is provable 
in the Hilbert-type $\mathbf{L_1}$, where $A_i$ $(1 \leq i \leq n)$ is  an atomic formula or 
a negated atomic one,

$(3)$ For any formula $A$ of $\mathbf{L_1}$, if $A$ is provable in the Hilbert-type $\mathbf{L_1}$, then 
it is provable in the the Gentzen-type $\mathbf{L_1}$ $($tableau method$)$, 

$(4)$ Contradiction and Dichotomy theorems for the Hilbert-type $\mathbf{L_1}$ hold,

$(5)$ No Hintikka formula is provable in the Hilbert-type $\mathbf{L_1}$, 

$(6)$ The Hilbert-type $\mathbf{L_1}$ is \L -decidable with respect to $\mathbf{HL_1}$ $($i.e. the set of 
all formula of $\mathbf{L_1}$ is the disjoint union of the set of all the theorem of the 
Hilbert-type $\mathbf{L_1}$ and that of all the theorem of $\mathbf{HL_1})$, where $\mathbf{HL_1}$ is 
the axiomatic rejection with Hintikka formulas as axioms.

Then we shall show that $\mathbf{HL_1}$ has the same strength with $\mathbf{HAR}$.

\rm
\end{abstract}

\smallskip

\noindent \small \it Keywords: \rm axiomatic rejection, the propositional fragment of Le\'{s}niewski's ontology, tableau method, refutation calculus, Hintikka formula, positive part, negative part, 
Sch\"{u}tte-type formulation of logic, cut elimination theorem, model, embedding, 
Russellian-type definite description. 
\normalsize


%

%

\section{Introduction}

In Ishimoto \cite{Ishimoto1977} and later Kobayashi-Ishimoto \cite{Kobayashi-Ishimoto1982} a logical system called  $\mathbf{L_1}$, was proposed as the 
propositional fragment of Le\'{s}niewski's ontoloty designated as  $\mathbf{L}$. The fragment $\mathbf{L_1}$ of 
Le\'{s}niewski's ontology is defined in its Hilbert-type version as the smallest class of formulas containg all the 
instances of tautology and the formulas of the form:

\smallskip 

1.1 $\enspace$ $\vdash \epsilon ab$ $\supset$ $\epsilon aa$,

\smallskip 

1.2  $\enspace$ $\vdash \epsilon ab$ $\wedge$ $\epsilon bc$. $\supset$ $\epsilon ac$,

\smallskip 

1.3b $\enspace$ $\vdash \epsilon ab$ $\wedge$ $\epsilon bb$. $\supset$ $\epsilon ba$,

\smallskip 

\noindent being closed under detachment as a rule (1.3 is due to Kanai \cite{Kanai1988}). 
Instead of 1.3b we may take, 
as an axiom, 

\smallskip 

1.3 (original) $\enspace$ $\vdash \epsilon ab$ $\wedge$ $\epsilon bc$. $\supset$ $\epsilon ba$,

\smallskip

\noindent which is the original one in Ishimoto \cite{Ishimoto1977}.

Le\'{s}niewski's (elementary) ontology  $\mathbf{L}$, on the other hand, is defined on the basis of the formula of 
the following form:

\smallskip 

1.4 $\enspace \enspace \vdash \epsilon ab \equiv. \exists x \epsilon xa \wedge \forall x \forall y (\epsilon xa 
\wedge \epsilon ya. \supset \epsilon xy) \wedge \forall x (\epsilon xa \supset \epsilon xb),$
\noindent or more simply, 

\smallskip 

1.5 $\enspace \enspace \vdash \epsilon ab \equiv. \exists x (\epsilon xa \wedge \epsilon xb) 
\wedge \forall x \forall y 
(\epsilon xa \wedge \epsilon ya. \supset \epsilon xy),$

\smallskip 

\noindent with first-order predicate logic (without equality) as underlying logic. 

The (well-formed) formulas of $\mathbf{L_1}$ to be referred to by such meta-logical variables 
as $A$, $B$, $\dots$ 
are defined in the well-known way in terms of $\epsilon$ (Le\'{s}niewski's epsion) and a (countable) infinite list of 
name variables, $a$, $b$, $\dots$ as well as a number of logical symbols sufficient for developing classical propositional 
logic and some auxiliary symbols. The (well-formed) formulas of $\mathbf{L}$ are defined analogously with quantifiers added. 

In what follows, all these symbols and their combinations will be employed only meta-logically. Outermost parentheses 
are always suppressed if no ambiguity arises therefrom.

As seen above we shall use a Polish-style notation such as $\epsilon ab$ for Le\'{s}niewski's epsilon as in Soboci\'{n}ski \cite{Sobocinski1967} 
 (for other notations, see a useful table in Simons \cite[p. 99]{Simons1987}).
 
 For Le\'{s}niewski's ontology in general refer, among others, to 
 Iwanu\'{s} \cite{Iwanus1973}, Lejewski \cite{Lejewski1958}, Luschei \cite{Luschei1962}, 
 Mi\'{e}ville \cite{Mieville1984}, Rickey \cite{Rickey1977}, Simons \cite{Simons1982, Simons1987}, 
 S\l upecki \cite{Slupecki1955}, Srzednicki-Rickey \cite{Srzednicki-Rickey1984}, 
 Surma et al. \cite{Surma-Srzednicki1992}, Stachniak \cite{Stachniak1981}, 
Cirulis \cite{Cirulis1977} and Urbaniak \cite{Urbaniak2014}. 
 
 Now the purpose of this paper is to prove, among others, that $\mathbf{L_1}$ is embedded in first-order predicate 
 logic with equality \it via \rm the translation $T$ to be specified below, that  is, embedding theorem:
 
\smallskip 
 
$\enspace \enspace  \enspace  \enspace \enspace \vdash_H \enspace A$

$\Longleftrightarrow 
\enspace TA \enspace \mbox{is valid in first-order predicate logic with equality, }$

\smallskip 

\noindent where $H$ is the Hilbert-type version of $\mathbf{L_1}$.

The translation $T$ which transforms every formula of  $\mathbf{L_1}$ into a formula of first-order 
 predicate logic with equality. The inductive definition $T$ is as follows. 
(This definition was first proposed by Ishimoto \cite{Ishimoto1977} on the basis of Prior \cite{Prior1965}.)

\smallskip 

T.1 $\enspace \enspace T \epsilon ab \enspace = \enspace F_b \iota xF_ax $,

\smallskip 

T.2  $\enspace \enspace T A \vee B \enspace  = \enspace T A \vee TB $,

\smallskip 

T.3 $\enspace \enspace T \sim  A  \enspace  = \enspace \sim T A $,

\smallskip 

\noindent $F_a, F_b, \dots$ are monadic predicate (variables) corresponding to name variables $a, b, \dots$ not 
necessarily, exhauting all of them. $F_b \iota xF_ax $, on the other hand, is the Russellian-type definite description and stands for:
$$\exists x(F_a x \wedge F_b x) \wedge \forall x  \forall y (F_a x \wedge F_a y. \supset x = y),$$
\noindent with the scope of the description confining to $F_b$. This embedding theorem was already proved in 
 Ishimoto \cite{Ishimoto1977} and later Kobayashi-Ishimoto \cite{Kobayashi-Ishimoto1982} by a different method. In the sequel, the theorem will be proved anew on the basis of a more general setting, which is summarized as the following (meta-)equivalences:
 
\smallskip 
 
$\vdash_T A \enspace  \Longleftrightarrow \enspace \vdash_H \enspace A$

$\enspace \enspace  \enspace  \enspace  \enspace \enspace  \Longleftrightarrow 
\enspace TA \enspace \mbox{is valid in first-order predicate logic with equality}$

$\enspace \enspace  \enspace  \enspace  \enspace \enspace  \Longleftrightarrow 
\enspace not \dashv_H A,$

\smallskip 

\noindent where $\vdash_T A \enspace (\vdash_H A$) signifies that $A$ is a thesis of the tableau method or the Gentzen-type version  (Hilbert-type version) of  $\mathbf{L_1}$, and $\dashv_H A$, 
 on the other hand, means that $A$ is axiomatically 
rejected in its Hilbert-type version: in other words, $A$ is an anti-thesis of it.  ($T A$ is the result from $A$ by applying the translation $T$ to $A$.  We will call it the \it $T$-transform of $A$. \rm  As easily seen, there are many formulas of predicate logic which are not a $T$-transform of a formula of $\mathbf{L_1}$. For other 
translation and embedding of $\mathbf{L_1}$, refer to Blass \cite{Blass1994}, Inou\'{e} 
\cite{Inoue1995b, Inoue2020, Inoue2020a}, Smirnov \cite{Smirnov1986, Smirnov1987} and Takano \cite{Takano1987}.

The paper to follow consists of nine sections with this introduction included with one appendix by 
the second authoer of this paper.

To begin with, the following second section will concern the tableau method version  of $\mathbf{L_1}$ 
and the proof of :
$$\vdash_T A \enspace \enspace  \Longrightarrow \enspace \enspace \vdash_H A,$$
along with an number of theorems which will turn out to be essential in the later development. In the third section, 
an axiomatic rejection $\mathbf{HAR}$ for $\mathbf{L_1}$, i.e. a logical system to derive all anti-theses of $\mathbf{L_1}$, 
will be introduced formally in its Hilbert-type version, whereas a Gentzen-type axiomatic rejection $\mathbf{GAR}$ for 
$\mathbf{L_1}$ will be proposed in 
the section seven. In the sections 4, 5, and 6, we shall prove on the basis of axiomatic rejection that: 
 
\smallskip 
 
$\enspace \enspace \enspace \enspace  \vdash_H \enspace A$

$\Longrightarrow 
\enspace TA \enspace \mbox{is valid in first-order predicate logic with equality}$

$\Longrightarrow 
\enspace not \dashv_H A.$

\smallskip 

For the first (meta-)implication, we can treat it as in Ishimoto \cite{Ishimoto1977}. For the second one, our argument for the proof is 
based on the model construction with assigments of finite subsets of the set of natural numbers to name variables as 
seen in Ishimoto \cite{Ishimoto1986}. (The original idea of the model construction is due to the third author of the present paper.)

It is emphasized that such a model for $\mathbf{L_1}$ is, without any change, one for $\mathbf{L}$ by making use of 
the arguments in Kobayashi-Ishimoto \cite{Kobayashi-Ishimoto1982} and that the model construced is finite. In the section six, the extension of a 
model for $\mathbf{L_1}$ to one for $\mathbf{L}$ will be used to prove a version of Separeation theorem (i.e. $\mathbf{L}$ 
is a conservative extention of $\mathbf{L_1}$), which was first proved in Ishimoto \cite{Ishimoto1977}. Our treatment about the model is a 
correction and a refinement of Ishimoto \cite{Ishimoto1986}.

Combining all those (meta-)implications with:
$$ not \dashv_H A \enspace \enspace \Longrightarrow \enspace \enspace \vdash_H A,$$
which will be proved in the third section (from Dichotomy theorem (Theorem 3.2)), we shall finally obtain the looked-for 
equivalences. The Dichotomy theorem and Contradiction theorem (Theorem 6.4 to be proved in the section six) for the Hilbert-type version of $\mathbf{L_1}$ provide us with the decidability of $\mathbf{L_1}$ 
(cf. S\l upecki \cite{Slupecki1972}).

In the seventh section we shall introduce a Gentzen-type axiomatic rejection for $\mathbf{L_1}$ and Contradiction theorem 
(Theorem 6.4) will be syntactically proved under the following postulate: no Hintikka formula of the form 
$A_1 \vee A_2 \vee \cdots \vee A_n \enspace (n \geq 1)$ is provable in the Hilbert-type version of $\mathbf{L_1}$, where 
$A_i \enspace (1 \leq i \leq n)$ is an atomic formula or a negated atomic one. The syntactical treatment of the theorem 
leads us to a novel syntactical proof of the cut elimination theorem for the tableau method version of $\mathbf{L_1}$, which 
will be carried out in the eighth section. The idea of such a proof would be applied to many logics. Our approach to prove 
cut elimination theorem was first explored in Inou\'{e}-Ishimoto \cite{Inoue-Isimoto1992} for classical propositional logic. 

Here, we wish to take an oppotunity of emphasizing that the proposed tableau method for $\mathbf{L_1}$ is a system to 
be developed within the bounds of its Hilbert-type version up until the sixth section. Such was the insight that S\l upecki and 
\L ukasiewicz had when they were working with the Aristotelian syllogistic in its Hilbert-type version. In fact, they saw through 
the Hilbert-type syllogistic a Gentzen structure hidden under the surface, although they never developed the structure as a 
self contained Gentzen-type logic. (This was attempted by 
Ishimoto-Kani-Kagiwada \cite{Ishimoto-Kanai-Kagiwada1984}, 
Kanai \cite{Kanai1990} for the Aristotelian syllogistic, and by 
Inou\'{e}-Ishimoto \cite{Inoue-Ishimoto-fc} for the Brentano-type syllogistic with Le\'{s}niewski's epsilon $\epsilon$.) Roughtly speaking, the Gentzen 
structure thus discovered was made use of very skillfully for the benefit of the Hilbert-type version of the logic concerned. 
And, the theses referring to the tableau mehtod should be understood only within the framework of $\mathbf{L_1}$ in its 
Hilbert-type version. Thus, $\vdash_T A$, for example, may be thought of, not only as a thesis of the Gentzen-type 
$\mathbf{L_1}$, but also as a theorem which belongs to the Hilbert-type $\mathbf{L_1}$ 
(see Theorem 2.2 in the following section). 

In the last section nine, as the chracterization theorem, we shall show: 

\noindent \textsc{Theorem} \rm (Characterization Theorem) \it 
 
\noindent The following statements are equivalent:

$(1)$ The Cut elimination theorem for the Gentzen-type $\mathbf{L_1}$ $($tableau method$)$ holds,

$(2)$  No Hintikka formula of the form $A_1 \vee A_2 \vee \cdots \vee A_n$ $(n \geq 1)$ is provable 
in the Hilbert-type $\mathbf{L_1}$, where $A_i$ $(1 \leq i \leq n)$ is  an atomic formula or 
a negated atomic one,

$(3)$ For any formula $A$ of $\mathbf{L_1}$, if $A$ is provable in the Hilbert-type $\mathbf{L_1}$, then 
it is provable in the the Gentzen-type $\mathbf{L_1}$ $($tableau method$)$, 

$(4)$ Contradiction and Dichotomy theorems for the Hilbert-type $\mathbf{L_1}$ hold,

$(5)$ No Hintikka formula is provable in the Hilbert-type $\mathbf{L_1}$, 

$(6)$  The Hilbert-type $\mathbf{L_1}$ is \L -decidable with respect to $\mathbf{HL_1}$ 
$($i.e. the set of all formula of 
$\mathbf{L_1}$ is the disjoint union of the set of all the theorem of the Hilbert-type $\mathbf{L_1}$ 
and that of all the theorem of $\mathbf{HL_1})$, where $\mathbf{HL_1}$ is 
the axiomatic rejection with Hintikka formulas as axioms. \rm
\smallskip

Then we shall show that the $\mathbf{HL_1}$ has the same strength with $\mathbf{HAR}$. $\mathbf{HL_1}$ was proposed in Inou\'{e} \cite{Inoue1995a}.

This paper contains an philosophical appendix by the second author, Arata Ishimoto. 

\section{Tableau method}

In this section, as in Kobayashi-Ishimoto \cite{Kobayashi-Ishimoto1982}, $\mathbf{L_1}$ will be developed by means of the tableau method, which in spite of 
its appearance may be understood in terms of the Hilbert-type $\mathbf{L_1}$ as remarked above. For this purpose, the notion of 
the positive and negative parts due to Sch\"{u}tte \cite{Schutte1960, Schutte1968, Schutte1977} will be introduced with a view to simplifying the subsequent development. 

\begin{definition} \bf 2.1 \rm 
The \it positive \rm and \it negative \rm parts of a formula $A$ are defined only as follows:

2.11 $\enspace A$ is a positive part of $A$,

2.12 $\enspace$ If $B \vee C$ ia a positive part of $A$, then $B$ and $C$ are positive parts of $A$,

2.13 $\enspace$ If $\sim B$ is a positive part of $A$, then $B$ is negaitive parts of $A$,

2.14 $\enspace$ If $\sim B$ is a negative part of $A$, then $B$ is positive parts of $A$.
\end{definition}

As suggested in this definition, the logical symbols to be employed in the sequel are $\vee$ (disjunction) and $\sim$ (negation) 
with other logical symbols being defined, if necessary, in their terms.

The specified occurrence of a formula $A$ as a positive (negative) part of another is indicated by $F[A_+]$ ($G[A_-]$) as 
exemplified below:

$F[A_+] = A$, 

$F[A_+] = \enspace \sim \sim A \vee B$, 

$F[A_+] = \enspace \sim B \vee (A \vee C)$, 

$G[A_-] = \enspace \sim A$,  

$G[A_-] = \enspace \sim \sim \sim A \vee B$,  

$G[A_-] = \enspace \sim \sim (\sim A \vee C) \vee B$, 

\noindent where all the formulas involved are assumed to  be different. Such expressions as $F[A_+, B_-]$ and the like are 
understood analogously with the proviso that specified formulas do not overlap with each other. 


On the basis of the above definition of the positive and negative parts of a formula, the \it tableaux \rm for $\mathbf{L_1}$ are 
defined by the following four reduction rules to be applied to a formula of $\mathbf{L_1}$:

$$\vee_{-} \enspace \frac{G [ A \vee B _{-} ]}{\enspace G [  A \vee B _{-} ] \vee \sim A \enspace | \enspace G [  A \vee B _{-} ] \vee \sim B \enspace }$$
\smallskip
$$\epsilon_{1} \enspace \frac{G [ \epsilon ab _{-} ]}{\enspace G [ \epsilon ab _{-} ] \vee \sim \epsilon aa \enspace }$$
\smallskip
$$\epsilon_{2} \enspace \frac{G [ \epsilon ab _{-}, \epsilon bc _{-} ]}{\enspace G [ \epsilon ab_{-}, \epsilon bc _{-} ] \vee \sim \epsilon ac \enspace }$$
\smallskip
$$\epsilon_{3b} \enspace \frac{G [ \epsilon ab _{-}, \epsilon bb _{-} ]}{\enspace G [ \epsilon ab_{-}, \epsilon bb _{-} ] \vee \sim \epsilon ba \enspace },$$
where all these reduction rules, as will be seen presently, should be understood as derived rules put up-side down  as far as we remain 
in the Hilbert-type $\mathbf{L_1}$. The rule $\epsilon_{3b}$ can be replaced by the following rule $\epsilon_{3}$ for the tableau mehtod:
$$\epsilon_{3} \enspace \frac{G [ \epsilon ab _{-}, \epsilon bc _{-} ]}{\enspace G [ \epsilon ab_{-}, \epsilon bc _{-} ] \vee \sim \epsilon ba \enspace }.$$

By reducing a formula by way of these reduction rules, a tableau is obtained for the formula. A branch of a tableau is \it closed \rm if it 
is ending with a formula of the form $F[A_+, A_-]$. A  tableau is said to be \it closed \rm if every branch of it is closed. A tableau is \it open \rm if it is not closed. A formula of $\mathbf{L_1}$ is \it provable \rm in the tableau 
method for $\mathbf{L_1}$ if there exists a closed tableau of it. It is also known that $\vee_{-}$ is sufficient 
for tableaux yielding all the instances of tautology. (For the formal definition of tableaux, consult 
Fitting \cite{Fitting1983} and Smullyan \cite{Smullyan1968}.)

The principal formulas of these rules, such as $A \vee B$ in the case of $\vee_{-}$, are minimal negative parts of 
the formulas to be reduced. Here, the \it minimal positive \rm or \it negaive parts \rm of a formula are the 
positive (negative ) parts of a formula which does not contain properly any positive or negative parts of the 
formula. The presence of the formulas repeated in the results of a reduction will be justified when we come across 
Hintikka formulas to be defined below.

For the purpose of illustration, the axioms 1.1, 1.2 and 1.3b of the Hilbert-type $\mathbf{L_1}$ will be 
proved by the proposed tableau method: 
$$\epsilon_{1} \enspace \frac{\sim \epsilon ab \vee \epsilon aa \enspace (= 1.1)}
{\enspace (\sim \epsilon ab \vee \epsilon aa) \vee \sim \epsilon aa \enspace ,}$$
\smallskip
$$\epsilon_{2} \enspace \frac{\sim \sim (\sim \epsilon ab \enspace \vee \sim \epsilon bc) 
\vee \epsilon ac \enspace (= 1.2)}
{\enspace (\sim \sim (\sim \epsilon ab \enspace \vee \sim \epsilon bc) \vee \epsilon ac) \vee \sim \epsilon ac 
\enspace \enspace , }$$
\smallskip
$$\epsilon_{3b} \enspace \frac{\sim \sim (\sim \epsilon ab \enspace \vee \sim \epsilon bb) 
\vee \epsilon ba \enspace (= 1.3b)}
{\enspace (\sim \sim (\sim \epsilon ab \enspace \vee \sim \epsilon bb) \vee \epsilon ba) \vee \sim \epsilon ba 
\enspace \enspace , }$$
where all these tableaux are closed. 


The following theorem (Theorem 2.3) fundamental in the subsequent development was stated and proved in 
Kobayashi-Ishimoto \cite{Kobayashi-Ishimoto1982}. But, here we take the opportunity of repeating it with a simple proof of it.

Before proceeding to the theorem, the notion of Hintikka formulas (for $\mathbf{L_1}$) is in order and it reads as 
follows: 

\begin{definition}\bf 2.2 \rm
A \it Hintikka formula \rm $A$ is a formula which satisfies the following conditions:

2.21 $\enspace$ $A$ is not of the form $F[B_+, B_-]$,

2.22 $\enspace$ If $A$ contains $B \vee C$ as a negative part of $A$, then 
it contains $B$ or $C$ as a negative part of it, 

2.23 $\enspace$ If $A$ contains $\epsilon ab$ as a negative part of $A$, then it contains $\epsilon aa$ 
as a negative part of it,

2.24  $\enspace$ If $A$ contains $\epsilon ab$ and  $\epsilon bc$ as a negative part of $A$, then it 
contains $\epsilon ac$ as a negative part of it,

2.25  $\enspace$ If $A$ contains $\epsilon ab$ and  $\epsilon bb$ as a negative part of $A$, then it 
contains $\epsilon ba$ as a negative part of it,
\end{definition}

We shall show some examples of Hintikka formula: 
$$\sim \epsilon ab \vee \epsilon ba \vee \sim \epsilon aa,$$
$$\sim (\epsilon ab \vee \epsilon bc)  \vee \sim \epsilon ab \vee \sim \epsilon aa,$$
$$\sim \epsilon ab \vee \sim \epsilon bc \vee \sim \epsilon ac \vee \sim \epsilon ba \vee 
\sim \epsilon aa \vee \sim \epsilon bb,$$
where some of the variables could be identical with each other. Here and in what follows, disjuncts 
are assumed to be associated in any way. 

\begin{theorem} \bf 2.1 \rm (Fundamental Theorem) \it 
Given a formula $($of $\mathbf{L_1})$, by reducing it by reduction rules there obtains a finite tableau, each 
branch of which ends either with a formula of the form $F[A_+, A_-]$ or with a Hintikka formula, whereby 
a branch is extended by a reduction rule only if the formula to be reduced is not of the form  $F[A_+, A_-]$ 
and the reduction gives rise to a formula not occurring in the formula to be reduced as negative part thereof. 
\end{theorem}

With a view to proving the Theorem, it is remarked in advance that there are only a finite number of 
subformulas of the given formula, and only some of them could be employed as a principal formula of a 
$\vee_{-}$ application. The \it principal formula \rm for $\vee_{-}$ is $A \vee B$, while that for $\epsilon_{1}$ 
($\epsilon_{2}$ and $\epsilon_{3b}$) is $\epsilon ab$ ($\epsilon ab$, $\epsilon bc$ and 
$\epsilon ab$, $\epsilon bb$).) There are also a finite number of the possible pairings of name variables in 
the given formula, and only some of them could be combined by the applications of  $\epsilon_{1}$, 
$\epsilon_{2}$ or $\epsilon_{3b}$ to yield a fresh occurrence of a negative part which was not in occurrence 
as such in the formula to be reduced. This proves the first half of the Fundamental Theorem. 

For proving the second half of the theorem, let us assume that extending a branch by way of reduction 
rules which is not ending with a formula of the form  $F[A_+, A_-]$, we come soon or later across a 
formula to which no rules are applicable any more without violating the requirement of the Theorem. 

We wish to show that the formula already constitues a Hintikka formula. If not, the formula would, for example, contain an $\epsilon ab$ as a negative part without containing another negative part $\epsilon aa$. We could, 
then, reduce the formula by $\epsilon _{1}$ against the assumption. The other properties of Hintikka formulas 
are taken care of analogously. 

From the requirement (in Theorem 2.1) for extending a branch, it immediately follows that any principal 
formula (or formulas) used as such before is (are) never employed again in the same status in the same 
branch.

A tableau, which is constructed in compliancewith the requirement, is said to be \it normal\rm . The tableaux 
obtained for the axioms of $\mathbf{L_1}$ are all normal if the name variables involved are different from 
each other.

In the sequel, we shall need the operation (due to Sch\"{u}tte) of removing a formula from another which 
contains the former as its positive of negative part. If a given formula is $F[A_+]$ ($G[A_{-}]$), the formula or the empty expression resulting by removing $A$ from the formula is denoted by $F[\enspace \enspace {  }_{+}]$ 
($G[\enspace \enspace {  }_{-}]$).
 
\begin{definition}\bf 2.3 \rm
Given $F[A_+]$ ($G[A_{-}]$), $F[\enspace \enspace {  }_{+}]$ ($G[\enspace \enspace {  }_{-}]$) is defined only 
as follows:

2.41 $\enspace$ If $F[A_+]$ is $A$, then $F[\enspace \enspace {  }_{+}]$ is the empty expression,

2.42 $\enspace$ If $F[A_+]$ is $F_1[A \vee B_+]$ or  $F_1[B \vee A_+]$, then $F[\enspace \enspace {  }_{+}]$ is
$F_1[B_{+}]$, 

2.43 $\enspace$ If $F[A_+]$ is $G_1[\sim A_-]$, then $F[\enspace \enspace {  }_{+}]$ is 
$G_1[\enspace \enspace {  }_{-}]$, 

2.44  $\enspace$ If $G[A_-]$ is $F_1[\sim A_+]$, then $G[\enspace \enspace {  }_{-}]$ is 
$F_1[\enspace \enspace {  }_{+}]$.
\end{definition}

The \it removal \rm of $A$ from $F[A_+]$ ($G[A_{-}]$) is defined by induction on the number of 
procedures 2.11-2.14 used for specifying $A$ as a positive (negative) part of $F[A_+]$ ($G[A_{-}]$). By the 
same induction, it follows that $F[\enspace \enspace {  }_{+}]$ ($G[\enspace \enspace {  }_{-}]$) is a 
well-formed formula or the empty expression, given $F[A_+]$ ($G[A_{-}]$).

The operation thus defined will be exemplified as follows.

If $F[A_+] = F_1[A \vee B_+] = F_2[(A \vee B) \vee C_+] = (A \vee B) \vee C$, then 
$F[\enspace \enspace {}_+] = F_1[B_+] = F_2[B \vee C_+] = B \vee C$.

If $F[A_+] = G_1[\sim A_-] = F_1[\sim \sim A_+] = F_2[\sim \sim A \vee B_+] = \sim \sim A \vee B$, then 
$F[\enspace \enspace {}_+] = G_1[\enspace \enspace {}_-] = F_1[\enspace \enspace {}_+] = 
F_2[B_+] = B$.

If $G[A_-] = F_1[\sim A_+] = G_1[\sim \sim A_-] = F_2[\sim \sim \sim A_+] = \sim \sim \sim A$, then 
$G[\enspace \enspace {}_-] = F_1[\enspace \enspace {}_+] = G_1[\enspace \enspace {}_-] = 
F_2[\enspace {}_+] = $ the empty expression.

If $G[A_-] = F_1[\sim A_+] = F_2[\sim A \vee B_+] = G_1[\sim (\sim A \vee B)_-] = $ 
$F_3[\sim \sim (\sim A \vee B)_+] = \sim \sim (\sim A \vee B)$, then 
$G[\enspace \enspace {}_-] = F_1[\enspace \enspace {}_+] = F_2[B_+] = $ $G_1[\sim B_-] = 
F_3[\sim \sim B_+] = \sim \sim B$.

\begin{lemma}\bf 2.1 \it $ $

\smallskip 
$\enspace \enspace \enspace \enspace \vdash_H \enspace F[A_+]  \equiv . F[\enspace \enspace {}_+] \vee A,  $

$\enspace \enspace \enspace \enspace \vdash_H \enspace G[A_-]  \equiv . G[\enspace \enspace {}_-] \vee \sim A. $
\end{lemma}

The Lemma is proved simaltaneously by induction on the number of procedures 2.11--2.14 applied for specifying $A$ as a positive (negative) part of $F[A_+]$ ($G[A_{-}]$). 

2.51 The basis is forthcoming right away, since we have $\vdash_H \enspace F[A_+]  \equiv . F[\enspace \enspace {}_+] \vee A$ with $F[\enspace \enspace {  }_{+}]$ being the empty expression, and the disjunction 
of the empty expression with any formula is identified with the formula, which is regarded as a stipulation. 

2.52 If $F[A_+]$ is $F_1[A \vee B_+]$ (or $F_1[B \vee A_+]$), then 
$\vdash_H F[A_+] \equiv F_1[A \vee B_+]$ (or $F_1[B \vee A_+]$) $\equiv. F_1[\enspace \enspace {  }_{+}] 
 \vee (A \vee B)$  (or $F_1[\enspace \enspace {  }_+] \vee (B \vee A)$) (by induction hypothesis) $. \equiv. (F_1[\enspace \enspace {  }_+] \vee A) \vee B 
. \equiv . F_1[B_{+}] \vee A$ (by induction hypothesis) 
$. \equiv . F[\enspace \enspace { }_{+}] \vee A$.
 
2.53  $\enspace$ If $F[A_+] = G[\sim A_-]$, then $\vdash_H F[A_+] \equiv G[\sim A_-] \equiv. 
G[\enspace \enspace {  }_{-}] \vee$ $\sim \sim A $ (by induction hypothesis) $.\equiv . 
G[\enspace \enspace {  }_{-}] \vee A .\equiv F[\enspace \enspace {  }_{+}] \vee A$.

2.54  $\enspace$ If $G[A_-] = F[\sim A_+]$, then $\vdash_H G[A_-] \equiv F[\sim A_+] \equiv. 
F[\enspace \enspace {  }_{+}] \vee$ $\sim A $ (by induction hypothesis) $.\equiv . 
F[\enspace \enspace {  }_{+}] \vee \sim A .\equiv G[\enspace \enspace {}_-] \vee \sim A$.
\begin{lemma}\bf 2.2 \it $ $

\smallskip 
$\enspace \enspace \enspace \enspace \vdash_H \enspace F[\enspace \enspace {}_+] \supset F[A_+],$

$\enspace \enspace \enspace \enspace \vdash_H \enspace G[\enspace \enspace {}_-] \supset G[A_-],$
\smallskip 

\noindent where $F[\enspace \enspace {  }_{+}]$ $(G[\enspace \enspace {  }_{-}])$ is not empty.
\end{lemma}
These two implications are proved on the basis of Lemma 2.1, respectively as follows:
$$\vdash_H \enspace F[\enspace \enspace {}_+] \supset .  F[\enspace \enspace {}_+]  \vee A , \enspace \vdash_H F[\enspace \enspace {}_+]  \vee A . 
\supset  F[A_+],$$
$$\vdash_H \enspace G[\enspace \enspace {}_-] \supset .  G[\enspace \enspace {}_-]  \vee \sim A , \enspace \vdash_H G[\enspace \enspace {}_-]  \vee \sim A . 
\supset   G[A_-].$$

The Lemma is to the effect of thinning in the sense of Gentzen [3]. 
Analogously, we have:
\begin{lemma}\bf 2.3 \it $ $

\smallskip 
$\enspace \enspace \enspace \enspace \vdash_H \enspace A \supset F[A_+],$

$\enspace \enspace \enspace \enspace \vdash_H \enspace \sim A \supset G[A_-].$
\smallskip
\end{lemma} 
\noindent The Lemma is again to the effect of thinning.
\begin{theorem}\bf 2.2 \it For any formula $A$ of $\mathbf{L_1}$, we have 
$$\vdash_T \enspace A \enspace \Longrightarrow \enspace \vdash_H \enspace A, $$
where $\vdash_T A $ means $($as already  stated$)$ that $A$ is a thesis of the tableau method version of 
$\mathbf{L_1}$ as interpreted in its Hilbert-type counterpart.
\end{theorem}
The theorem is proved by induction on the length of the tableau.

2.21 [Basis] The basis is taken care of by the following equivalence to be obtained on the basis of 
Lemma 2.1:
$$\vdash_H \enspace F[A_+, A_-] \equiv . \enspace F[\enspace \enspace {}_+, \enspace \enspace {}_-] \vee 
(A \vee \sim A),$$
the right-hand side of which is a tautology and, therefore, provable in the Hilbert-type version of 
$\mathbf{L_1}$.

2.22 [Induction steps] Induction steps are dealt with by the following equivalences all reduction rules, 
namely, $\vee_-$, $\epsilon_1$,  $\epsilon_2$ and  $\epsilon_{3b}$ (or $\epsilon_3$), 
which we are resorting, among others, to Lemmas 2.1 and 2.2 as well as to 1.1, 1.2 and 1.3b.

From $\enspace \vdash_H \enspace G[A \vee B_-] \enspace \equiv. \enspace G[\enspace \enspace {}_-] 
\vee \sim (A \vee B) \enspace$ and 

\noindent $\vdash_H \enspace G[\enspace \enspace {}_-] 
\vee \sim (A \vee B)\enspace \equiv. \enspace G[\enspace \enspace {}_-] \vee 
\sim (A \vee B) \vee \sim (A \vee B).$

$\enspace \enspace  \enspace \enspace \equiv. \enspace G[A \vee B_-] \vee \sim (A \vee B).$

$\enspace \enspace \enspace \enspace  \equiv . \enspace  
(G[A \vee B_-] \vee \sim A) \wedge  (G[A \vee B_-] \vee \sim B). \enspace  
(\mbox{for} \enspace \vee_-),$

\noindent we have 

$\vdash_H \enspace G[A \vee B_-] \equiv. 
(G[A \vee B_-] \vee \sim A) \wedge  (G[A \vee B_-] \vee \sim B). $
 
\smallskip

From $\enspace \vdash_H \enspace G[\epsilon ab_-] \vee \sim \epsilon aa. \enspace \supset . \enspace 
G[\epsilon ab_-]  \vee \sim \epsilon ab \enspace$ and 

$\vdash_H \enspace G[\epsilon ab_-] \vee \sim \epsilon ab . 
\enspace \equiv. \enspace G[\enspace \enspace {}_-] \vee 
\sim \epsilon ab \vee \sim \epsilon ab . $

$\enspace \enspace \enspace \enspace \equiv . \enspace  G[\epsilon ab_-] 
\enspace (\mbox{for} \enspace \epsilon_1),$

\noindent it follows that $\vdash_H \enspace G[\epsilon ab_-]  \vee \sim \epsilon aa. \enspace \supset \enspace  G[\epsilon ab_-].$

From 

\noindent $\enspace \vdash_H \enspace G[\epsilon ab_-, \epsilon bc_-]  \vee \sim \epsilon ac. \supset . 
G[\epsilon ab_-, \epsilon bc_-]  \vee \sim \epsilon ab \vee \sim \epsilon bc \enspace$ 

\noindent and 

$\vdash_H \enspace G[\epsilon ab_-, \epsilon bc_-]  \vee \sim \epsilon ab \vee \sim \epsilon bc.$

$\enspace \enspace \enspace \enspace \equiv. \enspace G[\enspace \enspace {}_-, \enspace \enspace {}_-] \vee 
\sim \epsilon ab \vee \sim \epsilon bc 
\vee \sim \epsilon ab \vee \sim \epsilon bc . $

$\enspace \enspace \enspace \enspace \equiv. \enspace G[\enspace \enspace {}_-, \enspace \enspace {}_-] 
\vee \sim \epsilon ab \vee \sim \epsilon bc \enspace (\mbox{for} \enspace \epsilon_2),$

\noindent we obtain $\enspace \vdash_H \enspace G[\epsilon ab_-, \epsilon bc_-]  \vee \sim \epsilon ac. \supset G[\epsilon ab_-, \epsilon bc_-].$
From 

\noindent $\enspace \vdash_H \enspace G[\epsilon ab_-, \epsilon bb_-]  \vee \sim \epsilon ba. \supset . 
G[\epsilon ab_-, \epsilon bc_-]  \vee \sim \epsilon ab \vee \sim \epsilon bc \enspace$ 

\noindent and 

$\vdash_H \enspace G[\epsilon ab_-, \epsilon bb_-]  \vee \sim \epsilon ab \vee \sim \epsilon bb.$

$\enspace \enspace \enspace \enspace \equiv. \enspace G[\enspace \enspace {}_-, \enspace \enspace {}_-] \vee 
\sim \epsilon ab \vee \sim \epsilon bb 
\vee \sim \epsilon ab \vee \sim \epsilon bb . $

$\enspace \enspace \enspace \enspace \equiv. \enspace G[\enspace \enspace {}_-, \enspace \enspace {}_-] 
\vee \sim \epsilon ab \vee \sim \epsilon bb \enspace (\mbox{for} \enspace \epsilon_3),$

\noindent we get $\enspace \vdash_H \enspace G[\epsilon ab_-, \epsilon bb_-]  \vee \sim \epsilon ba. \supset G[\epsilon ab_-, \epsilon bb_-].$

In view of the Theorem just proved, a proof in the tableau-method-type $\mathbf{L_1}$ is transformed 
into that of its correspondent in the Hilbert-type version of $\mathbf{L_1}$. In other words, a proof having 
tableau-method-proof is now thought of as a proof in the Hilbert-type $\mathbf{L_1}$ with each reduction 
rule to be understood as a derived rule of $\mathbf{L_1}$. 

The Theorem, it is remarked, constituetes the first step in obtaining the looked-for meta-equivalences as announced in the introduction. 
\begin{lemma}\bf 2.4 \it 
Every Hintikka formula contains at least one occurrence of atomic formulas either as its positive or 
negative part.
\end{lemma} 

Suppose, if possible, to the contrary. There would, then, be the shortest positive or negative part of the 
given Hintikka formula. If that be shorter positive parts against the assumption. If a formula were the 
shortest positive part of the form $A \vee B$, then $A$ and $B$ would the shorter positive parts 
against the assumption. If a formula having the shorter than $A \vee B$ again contrary to the hypothesis. 
If $\sim A$ were the shortest positive part, $A$ woul be a shorter negative part against the 
assumption. Lastly, if $\sim A$ were the shortest negative part, $A$ would be a shorter positive part again contrary to hypothesis. 

\section{Axiomatic rejection - its Hilbert-type version $\mathbf{HAR}$}

We are now in a position to state the axioms and rules for axiomatic rejection. 
The axiomatic rejection $\mathbf{HAR}$ to 
developed hereunder, it is noticed, is for the Hilbert-type version of $\mathbf{L_1}$. It constitutes a 
Hilbert-type axiomatic rejection in distinction to the Gentzen-type one to be introduced in what follows.

Axioms:

3.11 $\enspace \enspace \dashv_H \enspace \epsilon aa$, 

3.12 $\enspace \enspace \dashv_H \enspace \sim \epsilon aa$, 

\noindent where $a$ is a name variable specified for the purpose.

Rules:

3.13 $\enspace \enspace \vdash_H \enspace A \supset B, \enspace \dashv_H \enspace B 
\Longrightarrow \enspace \enspace \dashv_H \enspace A$, 

3.14 $\enspace \enspace \dashv_H \enspace A \enspace\Longrightarrow \enspace 
\enspace \dashv_H \enspace B$, 

\noindent in the second rule $A$ is obtained from $B$ by uniform substitution of a name variable for some 
occurring in $B$. (These two rules are due to \L ukasiewicz \cite{Lukasiewicz1951} 
and S\l upecki \cite{Slupecki1948, Slupecki1950}.)

3.15 $\enspace \enspace \dashv_H \enspace A \enspace\Longrightarrow \enspace 
\enspace \dashv_H \enspace A \vee \epsilon ab, \enspace \enspace $  (Kobayashi's rule)
\footnote{Kobayashi is the third author of the present paper.}

\noindent where $A$ is a Hintikka formula constituting a disjunction with all the disjuncts being 
either an atomic formula or a negatied atomic formula and $\epsilon ab$ does not occur in $A$ 
negated, i.e. as its negative part. 
Above $\dashv_H A$ means that $A$ is axiomatically rejected in the Hilbert-type 
$\mathbf{L_1}$, i.e. $\mathbf{HAR}$. Instead of it, we may denote it by 
$\vdash_{\mathbf{HAR}} A$. As easily seen, $A \vee \epsilon ab$ also constitues a Hintikka 
formula. (For axiomatic rejection, refer besides S\l upecki \cite{Slupecki1948, Slupecki1950}, 
\L ukasiewicz \cite{Lukasiewicz1951} and H\"{a}rtig \cite{Hartig1960}  
also to Goranko-Pulcini-Skura \cite{GPS2020}, Inou\'{e} \cite{Inoue1989, Inoue1991}, Ishimoto \cite{Ishimoto1981, Ishimoto1996}, 
Iwanu\'{s} \cite{Iwanus1973}, Skura \cite{Skura1990} and so on.)

We are now presenting an example of axiomatic rejection for the purpose of illustrating how our axioms 
and rules work in combination: 

3.2 $\enspace \enspace \dashv_H \enspace \epsilon ab \vee \epsilon bc. \supset \epsilon aa,$

\noindent where $a$, $b$ and $c$ are name variables different from each other.

(1) $\enspace \enspace \vdash_H \enspace \sim \epsilon bb \vee \sim \epsilon bb . 
\supset \sim \epsilon bb, \enspace \enspace $ tautology, 

(2) $\enspace \enspace \dashv_H \enspace \sim \epsilon bb \vee \sim \epsilon bb 
\enspace \enspace $ (1), 3.12, 3.13, 3.14, 

(3) $\enspace \enspace \dashv_H \enspace \sim \epsilon bc \vee \sim \epsilon bb 
\enspace \enspace $ (2), 3.14, 

(4) $\enspace \enspace \dashv_H \enspace (\sim \epsilon bc \vee \sim \epsilon bb) \vee \epsilon aa 
\enspace \enspace $ (3), 3.15, 

(5) $\enspace \enspace \vdash_H \enspace 
\sim (\epsilon ab \vee \epsilon bc) \vee \sim \epsilon bc \vee \epsilon aa. \supset . 
\sim \epsilon bc \vee \sim \epsilon bb \vee \epsilon aa$  

$\enspace \enspace $ $\enspace \enspace $  $\enspace \enspace $ tautology,

(6) $\enspace \enspace \dashv_H \enspace 
\sim (\epsilon ab \vee \epsilon bc) \vee \sim \epsilon bc \vee \epsilon aa 
\enspace \enspace $ (5), (4), 3.13, 

(7) $\enspace \enspace \vdash_H \enspace 
\sim (\epsilon ab \vee \epsilon bc) \vee \epsilon aa. \supset . 
\sim (\epsilon ab \vee \epsilon bc) \vee \sim \epsilon bc \vee \epsilon aa$  

$\enspace \enspace $ $\enspace \enspace $  $\enspace \enspace $ tautology,

(8) $\enspace \enspace \dashv_H \enspace 
\sim (\epsilon ab \vee \epsilon bc) \vee \epsilon aa 
\enspace \enspace $ (7), (6), 3.13.

\noindent From (8), we obtain 3.2.
\begin{lemma}\bf 3.1 \it For any formula $A$ of $\mathbf{L_1}$, we have:
$$\vdash_H \enspace A \equiv . \enspace A_1 \vee A_2 \vee \cdots \vee A_n 
\enspace \enspace (n \geq 1),$$
where $A_1, A_2, \cdots , A_n$ exhaust all the formulas which occur in $A$ as its minimal 
positive or negative parts with minimal negative parts prefixed with negation. $($The minimal positive 
$($negative$)$ 
parts of a formula, it is remembered, are those which contain neither positive nor negative parts of the 
formula except themselves.$)$
\end{lemma} 

This is proved by induction on the number of the minimal positive and negative parts which can be brought out by way of Lemma 2.1. 

The Lemma is exemplified as follows: 
$$\vdash_H \enspace 
\sim \sim \sim \epsilon ab \vee \sim \sim \epsilon ba \vee \sim \epsilon aa. \equiv . 
\sim \epsilon ab \vee \epsilon ba \vee \sim \epsilon aa,$$
$$\vdash_H \enspace \sim (\epsilon ab \vee \sim \epsilon bc) \vee 
\sim \sim \sim \epsilon ac \vee \sim \sim \epsilon bc \vee \sim \epsilon aa. $$
$$\equiv . \sim (\epsilon ab \vee \sim \epsilon bc) \vee 
\sim \epsilon ac \vee \epsilon bc \vee \sim \epsilon aa.$$

If $A$ is a Hintikka formula as is the case with the examples, at least, one formula among 
 $A_1, A_2, \cdots , A_n$ is atomic or the nagation of an atomic formula in view of Lemma 2.4. 
 
 We, next, wish to prove that every Hintikka formula of $\mathbf{L_1}$ is axiomatically 
rejected in $\mathbf{HAR}$, that is, a thesis of 
 $\mathbf{HAR}$. With this in view we are 
 proving a number of preparatory lemmas. 
\begin{lemma}\bf 3.3 \it Any formula of $\mathbf{L_1}$ of the form: 
$$\sim B_1 \vee \sim B_2 \vee \cdots \vee \sim B_n 
\enspace \enspace (n \geq 1),$$
is axiomatically rejected in $\mathbf{HAR}$, where $B_1, B_2, \cdots , B_n$ are atomic formulas. 
\end{lemma}

For proving the Lemma we substitute the name variable $a$ as specified in 3.12 for all the name variables 
occurring in the given formula. There, then, obtains,  
$$\sim \epsilon aa \vee \sim \epsilon aa \vee \cdots \vee \sim \epsilon aa, $$
which in turn is equivalent to  $\sim \epsilon aa$ (by classical propositional logic). From this it follows 
that the given formula is axiomatically rejected by 3.12, 3.13, 3.14. 
\begin{lemma}\bf 3.3 \it Any Hintikka formula of $\mathbf{L_1}$ of the form: 
$$A_1 \vee A_2 \vee \cdots \vee A_n \vee \sim B_1 \vee \sim B_2 
\vee \cdots \vee \sim B_m \enspace \enspace (n \geq 1, m \geq 1),$$
is axiomatically rejected in $\mathbf{HAR}$, where $A_1, A_2, \cdots ,$ $A_n, $ $B_1, B_2, \cdots , $ $B_m$ are atomic.
\end{lemma}

As shown by Lemma 3.2, 

3.31 $\enspace \enspace \sim B_1 \vee \sim B_2 \vee \cdots \vee \sim B_m \enspace \enspace (m \geq 1),$

\noindent is axiomatically rejected. The rule 3.15 for axiomatic rejection is, then, rejection of the given 
formula, whereby 3.31 is a Hintikka formula to begin with and the result of the application of 3.15 to a Hintikka formula again gives rise to another as remarked earlier in connection with the statement of the rule 3.15. 
\begin{lemma}\bf 3.4 \it Any Hintikka formula of $\mathbf{L_1}$ of the form: 
$$A_1 \vee A_2 \vee \cdots \vee A_n \enspace \enspace (n \geq 1),$$
is axiomatically rejected in $\mathbf{HAR}$, where $A_1, A_2, \cdots , A_n$ are atomic formulas. 
\end{lemma}

By substituting the variable $a$ as specified in 3.11 for all the variables taking place in the given formula, 
there obtains, 
$$\epsilon aa \vee \epsilon aa \vee \cdots \vee \epsilon aa, $$
which is equivalent to $\epsilon aa$ (by classical propositional logic). The given formula is again 
axiomatically rejected in view of 3.11, 3.13 and 3.14.

We are now in a position to treat the case of general Hintikka formulas. The next is one of our 
main results.

\begin{theorem}\bf 3.1 \rm (Basic Theorem) \it 

\noindent Every Hintikka formula of $\mathbf{L_1}$ is axiomatically rejected in $\mathbf{HAR}$.
\end{theorem}

For its proof, we shall begin with a given Hintikka formula is of the form 
as stipulated in the right-hand side of Lemma 3.1, namely,

\smallskip  

3.41 $\enspace \enspace A_1 \vee A_2 \vee \cdots \vee A_n 
\enspace \enspace (n \geq 1),$

\smallskip 

\noindent where, at least, one $A_i$ ($1\leq  i \leq n$) is atomic or the negation of an atomic formula. A Hintikka formula is always transformed into a Hintikka formula of the form 3.41 by Lemma 3.1 in view of Lemma 2.4 and the definition of 
Hintikka formula. 

The proof of Theorem 3.1 is carried out by induction on the number of the $A_i$s having the form $\sim (B_1 \vee B_2)$. 

The basis is to the effect that in 3.41 there does not occur any formula of the form $\sim (B_1 \vee B_2)$ 
and $\sim B_i$ $(i = 1$ or $i = 2)$. The basis holds from Lemmas 3.2, 3.3 and 3.4. Since 3.41 is a Hintikka formula, the given Hintikka formula 3.41 is equivalent to a formula of the form:
$$\cdots \vee \sim (B_1 \vee B_2) \vee \sim B_i \vee \cdots .$$

Now, by classical propositional logic, we have: 

\smallskip 

3.42 $\enspace \enspace \vdash_H \enspace (\cdots \vee \sim (B_1 \vee B_2) \vee \sim B_i \vee \cdots ) 
\supset (\cdots \vee \sim B_i \vee \cdots),$

\smallskip 

\noindent the consequence of which is again a Hintikka formula. 

By induction hypothesis we have: 

\smallskip 

3.43 $\enspace \enspace \dashv_H \enspace \cdots \vee \sim B_i \vee \cdots .$ 

\smallskip 

\noindent which in turn gives rise to, 

\smallskip 

3.44 $\enspace \enspace \dashv_H \enspace \cdots \vee \sim (B_1 \vee B_2) \vee \sim B_i \vee \cdots ,$ 

\smallskip

\noindent by 3.42, 3.43 and 3.13 as requested.

Since any Hintikka formula is equivalent to a formula of the form 3.41 by Lemmas 2.4 and 3.1, this 
completes the proof of our Basic theorem.

\begin{lemma}\bf 3.5 \it 
Given a branch of a tableau, which is ending with a Hintikka formula, every constituent formula of 
the branch is a positive part of the Hintikka formula, and such a formula implies the succedent one and is axiomatically rejected in $\mathbf{HAR}$.
\end{lemma} 

The first part is proved by induction on the length of the branch. The second and third parts are taken care 
of on the basis of the following these of Hilbert-type version $\mathbf{L_1}$ corresponding to reduction 
rules and the rule 3.13. 

$\vdash_H \enspace G[A \vee B_-] \enspace \supset . \enspace G[A \vee B_-] \vee \sim A, $

$\vdash_H \enspace G[A \vee B_-] \enspace \supset . \enspace G[A \vee B_-] \vee \sim B, $

$\vdash_H \enspace G[\epsilon ab_-] \enspace \supset . \enspace G[\epsilon ab_-] \vee \sim \epsilon aa, $

$\vdash_H \enspace  G[\epsilon ab_-, \epsilon bc_-] \enspace \supset . 
\enspace  G[\epsilon ab_-, \epsilon bc_-] \vee \sim \epsilon ac, $

$\vdash_H \enspace  G[\epsilon ab_-, \epsilon bb_-] \enspace \supset . 
\enspace  G[\epsilon ab_-, \epsilon bb_-] \vee \sim \epsilon ba, $

\noindent which were already mentioned in 2.22. 

\begin{corollary}\bf 3.1 \it 
Every positive $($the negation of negative$)$ part of a Hintikka formula of $\mathbf{L_1}$ is 
axiomatically rejected in $\mathbf{HAR}$.
\end{corollary}

This follows from Lemma 2.3 and Lemma 3.5 (or Theorem 3.1) and the rule 3.13.

By Theorem 2.1 (Fundamental Theorem) and Lemma 3.5, we have:
\begin{corollary}\bf 3.2 \it
Every formula of $\mathbf{L_1}$, which is not provable by the tableau method  for $\mathbf{L_1}$, is axiomatically rejected in $\mathbf{HAR}$.
\end{corollary}
\noindent From the Corollary and Theorem 2.2 there follows immediately an important theorem.
\begin{theorem}\bf 3.2 \rm $($Dichotomy theorem$)$ \it Every formula of $\mathbf{L_1}$, which is not a thesis of the Hilbert-type $\mathbf{L_1}$, 
is axiomatically rejected in $\mathbf{HAR}$ $($in its Hilbert-type version$)$.
\end{theorem}

In other words, every formula is provable or axiomatically rejected.

Suppose $A$ is not provable in the Hilbert-type version of $\mathbf{L_1}$. By Theorem 2.2, $A$ is also not 
provable in the tableau method for the Hilbert-type $\mathbf{L_1}$. By Corollary 3.2, $A$ is, then, 
axiomatically rejected in its Hilbert-type version.

Since the set of provable formulas as well as that of axiomatically rejected formulas in the Hilbert-type 
$\mathbf{L_1}$ are both recursively enumerable, Theorem 3.2, namely, Dichotomy theorem and 
Contradiction theorem (i.e. Theorem 6.4) to be proved in the sequel provides us with the 
decidability for the Hilbert-type version of $\mathbf{L_1}$.  Corollary 3.2 gives a decision procedure 
for the Hilbert-type version of $\mathbf{L_1}$. On the basis of a setting similar to ours, 
S\l upecki \cite{Slupecki1948, Slupecki1950} 
and \L ukasiewicz \cite{Lukasiewicz1951} gave a decision method for the Aristotelian syllogistic 
(cf S\l upecki \cite{Slupecki1972}). As will 
be seen in what follows, a Gentzen-type axiomatic rejection will give us a much more simpler decision 
method for $\mathbf{L_1}$.

Now, we are in a position to give a \it normal form \rm to each axiomatic rejection of the formula 3.2, 
which was already carried out by a more roundabout way.

With this in view, a (normal) tableau will be constructed for 3.2 as follows with all the name variables 
involved being different from each other: 
$$\displaystyle{\frac{\enspace \sim (\epsilon ab \vee \epsilon bc) \vee \epsilon aa \enspace \enspace 
( = 3.2)}
{
\displaystyle{\frac{\enspace  3.21 \enspace  }
{\enspace 3.22 \enspace } \enspace  (\epsilon_1) } \enspace  
 | \enspace  
\displaystyle{\frac{\enspace 
\sim (\epsilon ab \vee \epsilon bc) \vee \epsilon aa \vee \sim \epsilon bc \enspace
 (= 3.23) \enspace }
{\enspace \enspace 
\sim (\epsilon ab \vee \epsilon bc) \vee \epsilon aa \vee \sim \epsilon bc \vee \epsilon bb \enspace 
(= 3.24) \enspace } 
}\enspace  (\epsilon_1)}  (\vee_-)}$$
where 

\smallskip

$3.21 = \enspace \sim (\epsilon ab \vee \epsilon bc) \vee \epsilon aa \vee \sim \epsilon ab,$

$3.22 = \enspace \sim (\epsilon ab \vee \epsilon bc) \vee \epsilon aa  \vee \sim \epsilon ab \vee 
\sim \epsilon aa.$

\smallskip

We have by Lemma 3.5:
$$\vdash_H \enspace 3.2 \enspace \supset \enspace 3.23,$$
$$\vdash_H \enspace 3.23 \enspace \supset \enspace 3.24.$$

By Theorem 3.2, 3.24 is axiomatically rejected, since it is a Hintikka formula. In view of the above two 
implications and the rule 3.13, we are given an axiomatic rejection of 3.2 in the normal form, which 
consists in first obtaining a Hintikka formula and reject the given formula on its basis.

The above example for the normal form theorem can be regarded as a prototype for the argument 
developed in Inou\'{e} \cite{Inoue1989, Inoue1991}, in which axiomatic rejection with Hintikka formulas 
as axioms is proposed for classical propositional logic and its extensions. We will again touch on this in 
the section seven where a Gentzen-type axiomatic rejection is introduced for $\mathbf{L_1}$. 

\section{Translation and soundness}
In this section, we shall first recall the translation $T$ explained in the introduction. $T$ transforms 
every formula of $\mathbf{L_1}$ into a formula of first-order predicate logic with equality. The inductive 
definition of $T$ is as follows. (This definition was first proposed by Ishimoto \cite{Ishimoto1977} on 
the basis of Prior \cite{Prior1965}.)
\smallskip 

4.11 $\enspace \enspace T \epsilon ab \enspace = \enspace F_b \iota xF_ax $,

\smallskip 

4.12  $\enspace \enspace T A \vee B \enspace  = \enspace T A \vee TB $,

\smallskip 

4.13 $\enspace \enspace T \sim  A  \enspace  = \enspace \sim T A $,

\smallskip

\noindent $F_a, F_b, \dots$ are monadic predicate (variables) corresponding to name variables $a, b, \dots$ not 
necessarily, exhauting all of them. $F_b \iota xF_ax $, on the other hand, is the Russellian-type definite description and stands for 
$$\exists y(F_a x \wedge F_b x) \wedge \forall x  \forall y (F_a x \wedge F_a y. \supset x = y)$$
with the scope of the description confining to $F_b$. As easily seen, there are some formulas of predicate logic which are not a $T$-transform of a formula of $\mathbf{L_1}$. 
\begin{theorem}\bf 4.1 \rm (Soundness theorem) \it 
If $\vdash_H A$, then $TA$ is a thesis of first-order predicate logic with equality. 
\end{theorem} 

The proof is carried out by induction on the length of the proof in the HIlbert-type version of $\mathbf{L_1}$. 

The basis is taken care of on the basis of the following theses of predicate logic, respectively, 
corresponding to 1.11, 1.12 and 1.13:

\smallskip

$\vdash \enspace  F_b \iota xF_ax  \supset  F_a \iota xF_ax $,

$\vdash \enspace  F_b \iota xF_ax \wedge  F_c \iota xF_bx . \supset  F_c \iota xF_ax $,

$\vdash \enspace  F_b \iota xF_ax \wedge  F_b \iota xF_bx . \supset  F_a \iota xF_bx $,

\smallskip 

For treating the induction steps, let us assume that $TA$ and $TA \supset B$ ($TA \supset TB$) are 
provable in predicate logic. $TB$ is, then, forthcoming as a thesis of the logic by detachment.

This complete the proof of the Soundness theorem for the Hilbert-type $\mathbf{L_1}$. 

\section{Models for $\mathbf{L_1}$}
For proving, 
$$TA \enspace\mbox{is valid in first-order predicate logic with equality}$$
$\enspace \enspace  \Longrightarrow \enspace \enspace not \dashv_H A, $

\smallskip

\noindent which will be demonstrated in the next section as Theorem 6.1, we need some preparatory 
lemmas and definitions concerning the construction of models for $\mathbf{L_1}$, which are 
defined on the basis of the models for first-order predicate logic with equality. 
\begin{theorem}\bf 5.1 \it
For every Hintikka formula $A$, there is a model for $\mathbf{L_1}$ which falsifies $A$ and every 
positive (negative) part of it is false (true) there.
\end{theorem}

This is proved by defining a model for $\mathbf{L_1}$, where every atomic positive (negative) part of 
the Hintikka formula is made false (true), the existence of which is guaranteed by Lemma 2.4. 
The second part is taken care of by Lemma 2.3. But, we are proposing a method, which, though 
seemingly more complicated, will be convenient for obtaining a model for predicate logic for 
falsifying $TA$. The model in turn gives rise to the one for $\mathbf{L_1}$, as will be seen in what follows.
\begin{definition}\bf 5.1 \rm
A \it chain \rm a Hintikka formula $A$ is a (finite) collection of name variables $a_1, a_2, \dots , a_n$ 
($n \geq 1$) such that 

5.21 Every pair $a_i$ and $a_j$ ($1 \leq i \leq n, 1 \leq j \leq n$) belonging to the collection are 
connected by the relation defined as $\epsilon a_i a_j$ and $\epsilon a_j a_i$ both of which 
constitute negative parts of $A$,

5.22 The collection is maximal with respect to this property 5.21.
\end{definition}

As easily seen, the relation defined by 5.21 is reflexive, symmetric and transitive. 
\begin{definition}\bf 5.2 \rm
A \it tail \rm of a chain (of a Hintikka formula) is a name variable $b$ such that $\epsilon ab$ is a 
negative part of the Hintikka formula with the $a$, but not the $b$, being a member of the chain.
\end{definition}

For illustrative purposes, a number of Hintikka formulas will be presented with chains and tails associated 
thereto. The name variables involved, it is assumed, are different from each other.

5.31 $\enspace \enspace \sim \epsilon ab \vee \sim \epsilon ba \vee \sim \epsilon aa \vee \sim \epsilon bb,$

\noindent where $\{a, b\}$ is a chain without tails. ($\{a, b\}$, for example, is a set consisting of $a$ and $b$.) 

We shall here introduce a convenient notation for chains and tails. An expression $[x_1, x_2, \dots , x_n]$ 
means a chain $\{x_1, x_2, \dots , x_n\}$ of a Hintikka formula, where $x_1, x_2, \dots , x_n$ are name 
variables. By 
$$\begin{tikzpicture}[line width=0pt,node distance=2cm]
\node(a) []{$[x_1, x_2, \dots , x_n]$};
\node[right of=a,node distance=2.82cm](d) []{$y,$};
\draw (a) -- node[pos=.7, above]{} (d);
\end{tikzpicture}$$
we mean that $y$ is a tail of a chain $[x_1, x_2, \dots , x_n]$. With this notation, 5.31 is of the type $[a, b]$.

5.32 $\enspace \enspace \sim \epsilon ab \vee \sim \epsilon bc \vee \sim \epsilon ac \vee \sim \epsilon ba 
\sim \epsilon aa \vee \sim \epsilon bb,$

\noindent where $\{a, b\}$ is a chain with $c$ being its tail: 5.32 is of the type:
$$\begin{tikzpicture}[line width=0pt,node distance=2cm]
\node(a) []{$[a, b]$};
\node[right of=a,node distance=2.82cm](d) []{$c.$};
\draw (a) -- node[pos=.7, above]{} (d);
\end{tikzpicture}$$
5.33 $\enspace \enspace \sim \epsilon ab \vee \sim \epsilon bc \vee \epsilon ab \vee \sim \epsilon aa 
\vee \sim \epsilon bb,$

\noindent where $\{a\}$ and $\{b\}$, respectively, are different chains with c being $a$ tail common to 
these two chains:
$$
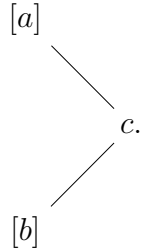
\begin{tikzpicture}[line width=0pt,node distance=2cm]
\node(a) []{$$};
\node(b) [above right of=a][]{$[a]$};
\node(c) [below right of=a][]{$[b]$};
\node[right of=a,node distance=2.82cm](d) []{$c.$};
\draw (b) -- node[midway, above]{} (d);
\draw (c) -- node[midway, above]{} (d);
\end{tikzpicture}$$
5.34 $\enspace \enspace \sim (\epsilon aa \vee \sim \epsilon bb) \vee \sim \epsilon ab 
\vee \sim \epsilon dc \vee \epsilon cb \vee \sim \epsilon bb \vee$

$\enspace \enspace \sim \epsilon ba \vee \sim \epsilon aa \vee \sim \epsilon dd \vee \sim \epsilon ae \vee \sim \epsilon be, $

\noindent where $\{a, b\}$ constitutes a chain with $e$ as its tail, while $\{d\}$ is another chain whose tail is $c$: 5.34 is of the type:

$\enspace \enspace \enspace \enspace \enspace \enspace \enspace \enspace  \begin{tikzpicture}[line width=0pt,node distance=2cm]
\node(a) []{$[a, b]$};
\node[right of=a,node distance=2.82cm](d) []{$e,$};
\draw (a) -- node[pos=.7, above]{} (d);
\end{tikzpicture}$

$\enspace \enspace \enspace \enspace \enspace \enspace \enspace \enspace \begin{tikzpicture}[line width=0pt,node distance=2cm]
\node(a) []{$[d]$};
\node[right of=a,node distance=2.82cm](d) []{$c.$};
\draw (a) -- node[pos=.7, above]{} (d);
\end{tikzpicture}$

5.35 $\enspace \enspace \sim \epsilon ab \vee \sim \epsilon ac \vee \epsilon bc 
\vee \epsilon da \vee \epsilon aa,$

\noindent where $\{a\}$ is the only chain having both $b$ and $c$ as its tails: 5.35  is of the type:
$$
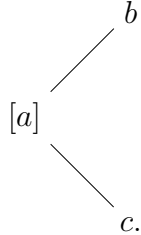
\begin{tikzpicture}[line width=0pt,node distance=2cm]
\node(a) []{$[a]$};
\node(b) [above right of=a][]{$b$};
\node(c) [below right of=a][]{$c.$};
\draw (a) -- node[midway, above]{} (b);
\draw (a) -- node[midway, above]{} (c);
\end{tikzpicture}$$

We are now going to describe a method of defining a model for $\mathbf{L_1}$ which falsifies the given 
Hintikka formula. We shall first confine ourselves to the definition of the models specific to the Hintikka 
formulas as above presented. All these models for $\mathbf{L_1}$ falsifies the given Hintikka formula, since
 every atomic formula taking place there as a positive (negative) part is false (ture) in these models, and 
 this makes every positive (negative) part of the Hintikka formula false (ture) as easily proved by induction 
 on the length of positive and negative parts (cf. Sch\"{u}tte \cite[p. 12 Theorem 1.6]{Schutte1977}). 
 
5.41 To begin with, we wish to define a model  $\mathcal{M}$ (for $\mathbf{L_1}$) which falsifies the 
Hintikka formula 5.31. 
 
The model $\mathcal{M}$ consists of two elements, namely, $\{1\}$ assigned to the members $a$ and $b$ 
of the chain and $\emptyset$ (the empty set) assigned to the infinite list of the remaining name variables. 
The truth value of atomic formula $\epsilon ab$ is that of $T\epsilon ab$, i.e. $F_b \iota F_ax$ to be 
defined on the basis of the model $\mathcal{M}'$ for first-order predicate logic with equality with the domain 
cconsisting$1$ only and equality standing for the identity between numbers. On the basis of this model, 
$\epsilon \{1\}\{1\}$ is true since $T\epsilon \{1\}\{1\}$ $=$ $\{1\} \iota x \{1\}x$ is true in $\mathcal{M}'$,  
while $\epsilon \{1\}\emptyset$, $\epsilon \emptyset \{1\}$ and $\epsilon \emptyset \emptyset$ are all 
seen false in $\mathcal{M}$, because their $T$-transforms are all false in $\mathcal{M'}$. The truth values of 
other formulas are defined on the basis of those for atomic formulas. (Here, $\{a\}$, for example, denotes a (monadic)predicate which is true only fo $a$, and $\emptyset$ the predicate constantly galse for any 
argument with respect to the model for the predicate logic. $\{a, b\}$, $\{a, b, c\}$ and the like are 
understood analogously.)

5.42 For defining a model (for $\mathbf{L_1}$) which falsifies the Hintikka formula 5.32, we assign $\{1\}$ to 
the members of the chain, el.e. $a$ and $b$, while the tail is given $\{1, n\}$ as its vvalue with $n$ being 
any natural number other that $1$, say $2$. The remaining name variables are assigned $\emptyset$ as before.

A model $\mathcal{M}$ (for $\mathbf{L_1}$ is, then, defined with $\emptyset$, $\{1\}$, $\{1, 2\}$ as the elements 
of the domain (the universe), while the truth value of any atomic formula is identified with that of its 
$T$-transform to be defined on the basis of the model $\mathcal{M}'$ for first-order predicate logic with 
$\{1, 2\}$ as its domain.

For example, the truth value of $\epsilon \{1\}\{1, 2\}$ is that of the following:
$$\exists x(\{1\}x \wedge \{1, 2\}x) \wedge \forall x \forall y (\{1\}x \wedge \{1\}y. \supset x = y),$$
whose truth vale is obtained by considering the following formula:

\smallskip

$(\{1\}1 \wedge \{1, 2\}1. \vee . \{1\}2 \wedge \{1, 2\}2)$

$\enspace \enspace \enspace \wedge \enspace (\{1\}1 \wedge \{1, 2\}1. \supset 1 = 1)$

$\enspace \enspace \enspace \wedge \enspace (\{1\}1 \wedge \{1, 2\}2. \supset 1 = 2)$

$\enspace \enspace \enspace \wedge \enspace (\{1\}2 \wedge \{1, 2\}1. \supset 2 = 1)$

$\enspace \enspace \enspace \wedge \enspace (\{1\}2 \wedge \{1, 2\}2. \supset 2 = 2).$

\smallskip  

\noindent The truth value of $\epsilon \{1\}\{1, 2\}$ is of course true. 

5.43 The Hintikka formula 5.34 is taken care of by assigning $\{1\}$, $\{2\}$, $\{1\}$, $\{1, 2, 3\}$ and 
$\emptyset$, respectively, to $a$, $b$, $c$ and the remaining variables and defining a model 
$\mathcal{M}$ (for $\mathbf{L_1}$ with its domain consisting of $\{1, 2, 3\}$ and $\emptyset$. The turth 
value of atomic formulas are again defined on the basis of a model  $\mathcal{M}'$ for first-order 
predicate logic with equality constructed analogously to the preceding two case with the domain 
consisting of $a$, $2$, $3$. $\mathcal{M}$ is a model for $\mathbf{L_1}$, and falsifies 5.33, since $T 5.33$ 
is false in $\mathcal{M}'$. 

5.44 The treatment of the Hintikka formula 5.34 proceeds by assigning  $\{1\}$, $\{1\}$, $\{2\}$, $\{1, 3\}$ 
$\{2, 4\}$, respectively, to  $a$, $b$, $d$, $e$, $c$ with the remainig variables assigned $\emptyset$. The 
$\mathcal{M}$ and $\mathcal{M}'$ are defined analogously to the preceding cases with the domin of 
$\mathcal{M}'$ consisting of $\{1\}$, $\{2\}$, $\{3\}$ and $\{4\}$ and 5.34 is false in $\mathcal{M}$, since 
$T 5.34$ is false in $\mathcal{M}'$. 

5.45 The Hintikka formula 5.35 is taken care of by assigning $\{1\}$, $\{1, 2\}$, $\{1, 3\}$ and 
$\emptyset$, respectively, to $a$, $b$, $c$ and the remaining variables. Everything goes as before, and 
$T 5.35$ is false in $\mathcal{M}'$. 

The definitions of models $\mathcal{M}$ and $\mathcal{M}'$ for $\mathbf{L_1}$ and first-order predicate 
logic with equality is respectively generalized in the follwoing way, given a Hintikka formula:

5.51 Every member of a chain is assigned one and the same unit set (of a natural number) with a 
different unit set assigned to a member of different chains. 

5.52 To tail we assign a set of natural numbers $\{m_1, m_2, \dots , m_k, N\}$, where $\{m_i\}$ 
($1 \leq i \leq k$) is the unit set associated with a member of a chain which is ending with the 
tail and $N$ is a number never employed so far in defining the model. The natural number $k$ depending
on the tail should be maximal. To make sure of the given situation, it is illustrated in the notation in 
5.31 as follows: 

\smallskip

Type
$$
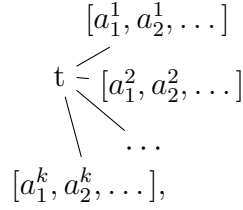
\begin{tikzpicture}
\node {t} 
[clockwise from=30, sibling angle=35]
child {node {$[a^1_1, a^1_2, \dots ]$}}
child {node {$[a^2_1, a^2_2, \dots ]$}}
child {node {$\cdots$}}
child {node {$[a^k_1, a^k_2, \dots ],$}};
\end{tikzpicture}$$

After the assignment
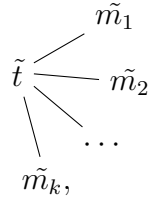
$$\begin{tikzpicture}
\node {$\tilde{t}$} 
[clockwise from=30, sibling angle=35]
child {node {$\tilde{m_1}$}}
child {node {$\tilde{m_2}$}}
child {node {$\cdots$}}
child {node {$\tilde{m_k},$}};
\end{tikzpicture}$$
where 

$\enspace \enspace \enspace \tilde{t} = \{m_1, m_2, \dots , m_k, N\}$,

$\enspace \enspace \enspace \tilde{m_1} = [\{m_1\}, \{m_1\}, \dots ]$,

$\enspace \enspace \enspace \tilde{m_2} = [\{m_2\}, \{m_2\}, \dots ]$,

$\enspace \enspace \enspace \tilde{m_k} = [\{m_k\}, \{m_k\}, \dots ]$.

\smallskip

5.53 To all other names, we assign the empty set $\emptyset$.

5.54 A model $\mathcal{M}'$ is, then, defined for first-order precicate logic with equality with the domain consisting of the (fine) set of natural numbers so far introduced. In case $\emptyset$ be the only 
set assigned to name variables, the domain of $\mathcal{M}'$ is any non-empty set of natural numbers. 
The truth value of any atomic formula, say $\epsilon ab$, in $\mathcal{M}$ is, then, identified with 
that of its $T$-transform $T\epsilon ab$, namely $F_b \iota x F_ax$ in $\mathcal{M}'$, and that of 
other formulas is defined on the basis of the truth values of atomic formulas. To the domain of $\mathcal{M}'$, 
any non-empty set of natural numbers could be adjoined without effecting the truth value in the model 
$\mathcal{M}$ (for $\mathbf{L_1}$). 

As expected, the model thus defined in general is the one for $\mathbf{L_1}$. With a view to proving this, 
all the entities obtained in the course of the model construction are classified into the following three 
categories, namely, 

$\emptyset$ $\enspace \enspace \enspace \enspace \enspace $ (the empty set as a predicate), 

$\{a\}$ $\enspace \enspace \enspace \enspace \enspace $ (a unit set as a predicate), 

$\{m_1, m_2, \dots , m_k, N\}$ $(k \geq 1)$ $\enspace \enspace \enspace $ 
(a finite set as a predicate), 

\noindent which are, respectively, assigned to name variables occurring neither as a members of 
a chain nor as a tail, members of a chain and tails.  

On the basis of such a classification, the $T$-transforms of 1.1, 1.2 and 1.3b turn out to be true 
in $\mathcal{M}'$, and the truth, then, gives rise to the satisfaction of axioms 1.1, 1.2 and 1.3b in 
$\mathcal{M}$ to be defined by way of $\mathcal{M}'$. 

5.61 To begin with, $T 1.1$ is seen to be true in $\mathcal{M}'$ for all the assignments to the name 
variables $a$ and $b$ as shown below: 

For $a = \emptyset$ or $\{m_1, m_2, \dots , m_k, N\}$, $T \epsilon ab$ is always false in $\mathcal{M}'$ 
irrespectively of any assignment to $b$ in view of the Russellian-type definite description. This makes 
$T 1.1$ and 1.1 true in $\mathcal{M}'$ and in $\mathcal{M}$, respectively. 

For $a = \{n\}$, $T \epsilon aa$ is true again by the definite description. This makes $T 1.1$ and, 
consequently, 1.1 true in $\mathcal{M}'$ and in $\mathcal{M}$, respectively. 

5.62 The axiom 1.2 is taken care of as follows: 

For $a = \emptyset$, $ \epsilon  ab$ is false in view of the definite description, and this gives rise to the 
truth of 1.2 in $\mathcal{M}$ through that of $T 1.2$. 

For $a = \{n\}$, and  $b = \emptyset$, $T 1.2$ is easily seen tru by the definite description.

For $a = \{n\}$, $b = \{m\}$ and $c = \emptyset$, $T 1.2$ is true in $\mathcal{M}'$ since $T \epsilon bc$ 
is false there making 1.2 true in $\mathcal{M}$. 

For $a = \{n\}$, $b = \{m\}$ and $c = \{1\}$, $T \epsilon ab$, $T \epsilon bc$ and $T \epsilon ab$ are 
all true if $n = m = 1$. This then, makes 1.2 true in $\mathcal{M}$. If $n =m$ and $m \neq 1$,  
$T \epsilon bc$ is false in $\mathcal{M}'$ making $T 1.2$ and 1.2 true in $\mathcal{M}'$ and in $\mathcal{M}$, respectively. For $n \neq m$, $T \epsilon ab$ is false and $T 1.2$ and 1.2 are true in $\mathcal{M}'$ and 
in $\mathcal{M}$, respectively. 

For $a = \{n\}$, $b = \{m\}$ and $c = \{m_1, m_2, \dots , m_k, N\}$, $T \epsilon ac$ is true if $n = m$ and 
$m$ is a member of $c$ , making $T 1.2$ and 1.2 true in $\mathcal{M}'$ and 
in $\mathcal{M}$, respectively. If $n = m$ and $m$ is not a member of $c$, then $T \epsilon bc$ is false, 
making  $T 1.2$ and 1.2 true in $\mathcal{M}'$ and in $\mathcal{M}$, respectively. If $n \neq m$, 
$T \epsilon bc$ is false and $T 1.2$ and 1.2 true in $\mathcal{M}'$ and in $\mathcal{M}$, respectively.

For $a = \{n\}$ and  $b = \{m_1, m_2, \dots , m_k, N\}$, $T \epsilon bc$ is false, making 
$T 1.2$ and 1.2 true in $\mathcal{M}'$ and in $\mathcal{M}$, respectively. 

For $a = \{m_1, m_2, \dots , m_k, N\}$, $T \epsilon bc$ is always false, for any value of $b$, and $T 1.2$ and consequently, 1.2 is ture in $\mathcal{M}'$ and in $\mathcal{M}$, respectively. 

5.63 The axiom 1.3 is seen to be true in $\mathcal{M}$ in the following way: 

For $b = \emptyset$ or $b = \{m_1, m_2, \dots , m_k, N\}$, $T \epsilon bb$ is false in $\mathcal{M}'$  
making $T 1.3$ nad 1.3 true, respectively, in $\mathcal{M}'$ and in $\mathcal{M}$.  

For $b = \{n\}$, $a = \emptyset$, $T \epsilon ab$ is false in $\mathcal{M}'$ making 1.3 true in $\mathcal{M}$. 

For $b = \{n\}$, $a = \{m\}$ and $n = m$, $T \epsilon ba$ is true in $\mathcal{M}'$ and this makes $T 1.3$ 
and 1.3 true respectively, in $\mathcal{M}'$ and in $\mathcal{M}$.

For $b = \{n\}$, $a = \{m\}$ and $n \neq m$, $T \epsilon ab$ is false in $\mathcal{M}'$. This makes $T 1.3$ 
and 1.3 true respectively, in $\mathcal{M}'$ and in $\mathcal{M}$. 

5.64 This completes the proof that the axioms 1.1--1.3 for $\mathbf{L_1}$ are satisfied by the model 
$\mathcal{M}$ since we have $T A \supset B \equiv . TA \supset TB$. Therefore, this model constructed 
by 5.51--5.54 constitutes a finite model for $\mathbf{L_1}$. 

5.65 We next have to check whether such a model (for $\mathbf{L_1}$) as constructed above actually 
falsifies the given Hintikka formula. By the definition of Hintikka formula, we may only consider atomic 
positive (negative) parts of it. In other words, if each atomic positive (negative) part of it is assigned 
falsity (truth) (via $\mathcal{M}'$), then the Hintikka formula sis false in $\mathcal{M}$ (via $\mathcal{M}'$) 
(cf, for example Sch\"{u}tte \cite[p. 12, Theorem 1.6]{Schutte1977}). (Note that such atomic positive (negative) parts of it 
are the minimal positive (negative) ones of it, but they do not always exhaust all the minimal ones in some 
cases.)

5.66 If the Hintikka formula has $T \epsilon ab$ as its positive part, the n we have the following 
possibilities as the result of our assignment: 

\smallskip

$\epsilon \emptyset \emptyset$, 

$\epsilon \emptyset \{n\}$, 

$\epsilon \emptyset \{m_1, m_2, \dots , m_k, N\}$, 

$\epsilon \{n\} \emptyset$, 

$\epsilon \{p\} \{q\}$, 

$\epsilon \{r\} \{m_1, m_2, \dots , m_k, N\}$, 

$\epsilon \{m_1, m_2, \dots , m_k, N\} \emptyset$, 

$\epsilon \{m_1, m_2, \dots , m_k, N\} \{n\}$, 

$\epsilon \{m_1, m_2, \dots , m_k, N\} \{m'_1, m'_2, \dots , m'_j, L\} $, 

\smallskip

\noindent which are all false in $\mathcal{M}$, where $p \neq q$ and $r$ is not a member of 
$\{m_1, m_2, \dots , m_k, N\}$. The other possibilities do not happen because of the definitions 
2.2, 5.2, 5.3, 5.51--5.53. 

5.67 If the Hintikka formula has $\epsilon ab$ as its negative part, then we have the only 
following possibilities as the result of our assignment: $\epsilon \{n\} \{n\}$ and 
$\epsilon \{m\} \{m_1, m_2, \dots , m_k, N\}$, which are true in $\mathcal{M}$, where $m$ is an element of 
$\{m_1, m_2, \dots , m_k, N\}$. By similar reasoning as in 5.66, the other possibilities do not happen. It is, 
however, remarked that the Hinkikka formula does not contain formulas of the form $\epsilon tb$ as its 
negative part, where $t$ is a tail and $a$ is an arbitrary name variable. If such a formula were a negative 
part of the Hintikka formula, then the Hintikka formula would be of the form $F[\epsilon ta_-, \epsilon tt_-]$ 
by the definition of Hintikka formula. Thus $t$ is an element of some chain. This contradicts the following.
\begin{proposition}\bf 5.1 \it
Suppose that $A$ is a Hintikka formula of $\mathbf{L_1}$. Then no tail of a chain of $A$ belongs 
to other chains of it.
\end{proposition} 

For proving the Proposition, let $C_1$ and $C_2$ be chains of $A$ and $b$ a tail of $C_1$. Suppose 
that $b$ is a member of $C_2$. Then there is a name variable of $C_1$ such that $\epsilon ab$ is 
a negative part of $A$. Since $b$ is an element of $C_2$, $A$ contains  $\epsilon bb$ as its negative 
part. By 2.25, $A$ thus contains  $\epsilon ba$ as its negative part. In other words, $b$ is a member of 
$C_1$, which contradicts the definition of tail.

5.68 From 5.68--5.67, it follows that the given Hintikka formula is falsified in $\mathcal{M}$ and that 
every positive (negative) part of it is false (true) in $\mathcal{M}$ as understood by induction of the length of positive (negative) part. This is the second part of Theorem 5.1. 

By 5.65--5.67, we complete the proof of Theorem 5.1. (The original idea of our model construction is due 
to the third author of the present paper.) 

We are, now, taking the opportunity of demonstrating that the model $\mathcal{M}$ as above defined for 
$\mathbf{L_1}$ also constitutes the one for $\mathbf{L}$ as well. This will play an important role, in 
particular, with reference to the proof of Separation theorem to be shown in the section six.

The model thus constructed for falsifying the given Hintikka formula is a model 
$\mathcal{M} =  <\mathcal{D}, \epsilon>$ such that $\mathcal{D}$ is a finite set of subsets of the set of 
natural numbers, which are regarded as Le\'{s}niewskian names, and the truth values of the formulas 
are reduced to those of their $T$-transforms, which are in turn based upon a model 
$\mathcal{M}' =  <\mathcal{D}', F_a, F_b, \dots>$ for first-order predicate logic with equality. Here, 
$\mathcal{D}'$ is non-empty set of natural numbers, not necessary infinite. Secondly and more importantly 
$F_a, F_b, \dots$ do not exhaust all the subsets of $\mathcal{D}'$. Put it the other way round, a number 
of subsets of $\mathcal{D}'$ are remaining anonymous with being named by any names.

Such a situation is sometimes responsible for the appearance of the so-called singular names 
which are taking place in the process of defining models for $\mathbf{L_1}$. A \it singular name\rm , which 
was introduced in Ishimoto \cite{Ishimoto1977}, is a name which is not an atom, but contains only one name in the sense of 
$\epsilon$-relation, where an atom \rm means a member $a$ of the domain of a model $\mathcal{M}$ 
such that $\epsilon aa$ is true in $\mathcal{M}$. Here, note that a name used in Ishimoto \cite{Ishimoto1977} is a subset 
of the set of natural numbers, i.e. an element of the domain of a model for $\mathbf{L_1}$ in our 
present contex. A singular name is occurring, for example, in the model as defined  in 5.42, where $b$, 
which happens to be a tail, is singular wihtout being an atom, but containing only one name, namely, $a$. 
As will be seen presently, such a model fails to be a model for $\mathbf{L}$. 

The remedy is not difficult to think of. It is only requested to adjoin $\{N\}$s with $N$s, each of which 
is contained in the counterpart $\{m_1, m_2, \dots , m_k, N\}$ of each given singular name. In the 
case of 5.42, in addition to $a$ and $b$ there takes place another name, say, $c$ which is the unit 
set $\{2\}$. 

The remedy skeched above, which makes use of introducing some new names, is based on the argument 
in Ishimoto \cite{Ishimoto1977} and later in Kobayashi-Ishimoto \cite{Kobayashi-Ishimoto1982}, which will be discussed 
below in our context. The reader will see that the argument can be applied to our case without any 
change, because our models constructed above coincide with ones defined by the model construction 
in Kobayashi-Ishimoto \cite{Kobayashi-Ishimoto1982}, if we identify names in the just cited paper with our subsets of the set of 
natural numbers: in other words, every atomic positive (negative) part of a given Hintikka formula is 
false (true) and the rest of all atomic formulas is falsified in both models. This is the truth concept for 
atomic formulas when we argue about the models for $\mathbf{L}$, while the counterpart of the models 
for $\mathbf{L_1}$ is on the definite description. This identification is the very trick of connecting the 
argument of ours to that of the paper cited above. (The idea is due to the second author of the 
present paper.) To continue the remedy, we shall show the following strightforward 
\begin{lemma}\bf 5.1 \it 
Given a model $(\mathbf{L_1})$, another model is defined which does not involve singular names, 
and the truth value of the formula in the original model are remaining the same in the new model.
\end{lemma}  

It is remarked in passing that any model thus constructed always contains atoms. This is because 
such a model involves names corresponding to unit sets, and names of unit sets are atoms. 

For demonstrating that the models for $\mathbf{L_1}$ thus augmented, if necessary, 
are the models for $\mathbf{L}$ as well, we need in advance the following preparatory theorem.
\begin{theorem}\bf 5.2 \it 
Given a model for $\mathbf{L_1}$, it is also the one for $\mathbf{L}$ if and only if it does not involve 
any singular names.
\end{theorem} 

Necessity: Given a model $\mathcal{M} =  <\mathcal{D}, \epsilon>$ for $\mathbf{L_1}$, suppose it 
contains a singular name $b$, of which the only element is a (in the sense of $\epsilon$-relation) 
and it is, of course, an atom.

Then, 

(1) $\enspace \enspace \exists x (\epsilon xb \wedge \epsilon xb)$, 

\noindent is ture in the model $\mathcal{M} $. This is because 

(2) $\enspace \enspace \epsilon ab \wedge \epsilon ab$, 

\noindent is true in $\mathcal{M}$ with $a$ and $b$ belonging to $\mathcal{D}$.

We also have the truth of 

(3) $\enspace \enspace \forall x \forall y (\epsilon xb \wedge \epsilon yb. \supset \epsilon xy)$, 

\noindent in $\mathcal{M}$, since $a$ is the only name ($\in \mathcal{D}$) such that $\epsilon ab$ 
is tre. $\epsilon bb$, on the other hand, is false, since $b$ is singular not being an atom. This makes the 
axiom schema 1.5 for $\mathbf{L}$ false for $a = b$ in the model $\mathcal{M}$.

Sufficiency: Given a model $\mathcal{M} =  <\mathcal{D}, \epsilon>$ for $\mathbf{L_1}$. Suppose 
$\mathcal{M}$ is not a model for $\mathbf{L}$, although it is for $\mathbf{L_1}$. Then, there are, at 
least, two names $a$ and $b$. This is to the effect that 

(4) $\enspace \enspace \epsilon ab$, 

\noindent is false with 

(5) $\enspace \enspace \exists x (\epsilon xa \wedge \epsilon xb) \wedge 
\forall x \forall y (\epsilon xa \wedge \epsilon ya. \supset \epsilon xy)$, 

\noindent being true in the model. In this connection, it is remarked, the converse implication 
$(4) \supset (5)$ is true in the model for $\mathbf{L_1}$.

Since (5) is true in the model, there is an $x$ ($\in \mathcal{D}$) such that 

(6) $\enspace \enspace \epsilon xa \wedge \epsilon xb$, 

\noindent is true there.

Further, suppose, if possible, that there is an $y$ ($\in \mathcal{D}$) such that 

(7)  $\enspace \enspace \epsilon ya \wedge \sim \epsilon xb$.

Inview of (6), (7) and (5), $\epsilon yx$ is true in the model, which gives rise to the truth of $\epsilon yb$ 
by 1.2. This, however, contradicts (7).

We, thus, have the truth in the model of 

(8) $\enspace \enspace \forall x  (\epsilon xa. \supset \epsilon ab)$, 

\noindent from which follows by (4) that 

(9) $\enspace \enspace \epsilon aa$, 

\noindent is not true in the model. This is to the effect that $a$ is a singular name of the model 
in view of (5).

This completes the proof of Theorem 5.2. 

From the Theorem there straightforwardly obtain:
\begin{theorem}\bf 5.3 \it 
A model constructed as above for $\mathbf{L_1}$, if properly extended when necessary, constitutes 
at the same time one for $\mathbf{L}$. 
\end{theorem}

In the above argument we followed the usual model-theoretic interpretation for Le\'{s}niewskian 
quantification, while an alternative interpretation, namely the substitutional one has been much discussed 
so far by some of the leading philosophers. In our setting we are assigning to each name variable a 
certain element of a structure, i.e. a subset of the set of natural numbers, whereas the substitutional 
interpretation does not assign anything to the name variables and thus does not need any domain of 
our model or the like. Because of the reason the substitutional interpretation would appeal to some, 
although it could not be the only reason. We will, here, not go into the alternative interpretation further. 
But, we shall cite the literature about the interpretation, where we see a variety of arguments or it, as 
follows: K\"{u}ng \cite{Kung1974, Kung1977}, K\"{u}ng-Canty \cite{Kung-Canty1977}, 
Quine \cite{Quine1969}, Rickey \cite{Rickey1985}, Simons \cite{Simons1985, Simons1985a} and so on.

As a concluding remark of the present section, it is emphasized that the domain of the model for 
$\mathbf{L_1}$ and $\mathbf{L}$  is a finite set whose elements are all finite sets and that our argument is treated within the bounds of first-order logic. In addition, we mention that the treatment for the remedy of 
the models with singular names is not the only one: we can, for example, make use of the result of 
Takano \cite{Takano1985} for an alternative treatment. (The idea for the alternative treatment is due to the first 
author of the present paper. We decide to take the second author's idea for the present paper.) 
Takano's paper contains a proof of the completeness theorem for Le\'{s}niewski's ontology 
$\mathbf{L}$ with 
respect to a natural truth concept, 
$$\epsilon ab \enspace \mbox{is true (in a structure)} 
\Longleftrightarrow \exists p(a =\{p\} \wedge p \in b),$$
where the right $\in$ means the membership relation of the set theory, which is similar to ours. The 
embedding of Le\'{s}niewski's ontoloty into the monadic second-order predicate logic in Smirnov [59] may be regarded as a syntactical version of such a natural interpretation of Le\'{s}niewski's epsilon (cf. 
Takano \cite{Takano1987}).

For a similar model construction for the Aristotelian syllogistic, one may consult Kanai \cite{Kanai1990}. 

In addition, we shall here mention some application of $\mathbf{L_1}$ for natural 
language as follows: 
Ishimoto \cite{Ishimoto1982, Ishimoto1986, Ishimoto1990}, Ishimoto-Shimidzu 
\cite{Ishimoto-Shimidzu1986} and Shimidzu \cite{Shimidzu1989, Shimidzu1990}. 

\section{Axiomatic rejection and embedding theorem} 
With a view to proving the theorem announced at the beginning of the preceding section, let us assume 
that $\dashv_H A$. We, then, wish to prove that $TA$ is not valid in first-order predicate logic with 
equality. 

The proof is carried out by induction on the number of rules applied for axiomatically rejecting $A$.

6.1 The basis does not present any difficulties, since both $\epsilon aa$ and $\sim \epsilon aa$ 
constitute Hintikka formulas, and their $T$-transforms are both falsified by a model for first order 
predicate logic with equality in view of Theorem 5.1. In this connection, it is remembered, the model 
for $\mathbf{L_1}$, which falsifies the Hintikka formula, is defined on the basis of a model for 
predicate logic, in which the $T$-transform of the Hintikka formula is also false.

6.2 Induction steps:  

6.21 The last applied rule for rejection is 3.13. By induction hypothesis, there is a model for first-order 
predicate logic which falsifies $TB$. This is to te effect that $TB$ is not  a thesis of first-order predicate 
logic by the soundness theorem for the logic. Now, we have $\vdash_H \enspace A \supset B$. By 
Soundness theorem (Theorem 4.1), $TA \supset B = TA \supset TB$ is a thesis of first-order 
predicate logic with equality. From this it follows that $TA$ is not valid in the predicate logic along with $TB$.

6.22 The last applied rule for rejection is 3.14. For taking care of this case, let us assume that $T$ is 
falsified by a model for first-order predicate logic with equality, and $A$ is obtained from $B$ by uniform 
substitution for some name variables occurring in $B$. As easily seen, $TB$ is also falsified by the 
same model by identifying the value of $b$ with that of $a$, where $a$ is substituted for $b$ in $B$.

6.23 The last applied rule for rejection is 3.15. The case is taken care of without resorting to 
induction hypothesis. In fact, not only the given Hintikka formula, but $A \vee \epsilon ab$ is also a Hintikka formula, and its $T$-transform is made false by a model for first-order predicate logic with equality by 
Theorem 5.1.

This completes the proof of:
\begin{theorem}\bf 6.1 \it 
$$TA \enspace \mbox{is valid in first-order predicate logic with equality}$$
$\enspace \enspace  \Longrightarrow \enspace \enspace not \dashv_H A.$
\end{theorem}
Making use of the completeness theorem of first-order predicate logic, this Theorem together with 
Theorems 2.2 and 4.1 as well as with Corollary 3.2 gives rise to the looked-for equivalences:
\begin{theorem}\bf 6.2 \rm (Main Theorem) \it  

$\vdash_T A \enspace  \Longleftrightarrow \enspace \vdash_H \enspace A$

$\enspace \enspace  \enspace  \enspace  \enspace \enspace  \Longleftrightarrow 
\enspace TA \enspace \mbox{is valid in first-order predicate logic with equality}$

$\enspace \enspace  \enspace  \enspace  \enspace \enspace  \Longleftrightarrow 
\enspace not \dashv_H A.$
\end{theorem} 
\noindent which were announced in the first section. 

In particular, from the second equivalence, there obtains a theorem to the effect that $\mathbf{L_1}$ is 
embedded in first-order predicate logic with equality \it via \rm the translation $T$, namely, 
\begin{theorem}\bf 6.3 \it 
$\vdash_H A$ if and only if $TA$ is a thesis of first-order predicate logic with equality.
\end{theorem}

Here, we are again making use of the completeness of the predicate logic.

In the third section, Dichotomy theorem (Theorem 3.2) was proved to the effect that every formula (of 
$\mathbf{L_1}$ is either provable or axiomatically rejected in the Hilbert-type $\mathbf{L_1}$. Now, another 
theorem to be coupled with this Theorem will be proved in this section as already mentioned in the 
section three. It will be called Contradiction theorem.
\begin{theorem}\bf 6.4 \rm (Contradiction theorem) \it It is not the case that for any formula $A$ $($of  $\mathbf{L_1})$, 
$\vdash_H A$ and $\dashv_H  A$ at the same time.
\end{theorem}

The proof is straightforward in view of the second and third equivalences of Theorem 6.2. A syntactical 
proof is also possible. I will be presented in the section eight (under some postulate).

Availing ourselves of Contradiction theorem just proved, there is forthcoming:
\begin{corollary}\bf 6.1 \it For any formula $A$ of $\mathbf{L_1}$, we have:

$\dashv_H  A$ if and only if $TA$ is not valid in first-order predicate logic with equality.
\end{corollary} 

Suppose $\dashv_H  A$, then we have not $\vdash_H  A$ by Contradiction theorem, which in turn 
gives rise to that $TA$ is not valid in the predicate logic by Theorem 6.3. If not $\dashv_H  A$, then we 
have $\vdash_H  A$ by Dichotomy theorem (Theorem 3.2), which, then yields the negation of the right 
side of the Corollary by Theorem 4.1 (Soundness theorem).

The following theorem is a version of Separation theorem, which was first proved in 
Ishimoto \cite[Theorem 3.4, p. 293]{Ishimoto1977}. 
\begin{theorem} \bf 6.5 \rm (Separation theorem) \it If a quantifier-free formula $A$ of $\mathbf{L}$, i.e. a formula 
belonging to $\mathbf{L_1}$ is valid, then $A$ is already a thesis of $\mathbf{L_1}$. 
\end{theorem}

In other words, $\mathbf{L}$ is a conservative extension of $\mathbf{L_1}$. 

Suppose, if possible, $A$ is not provable in the Hilbert-type $\mathbf{L_1}$. Then, by Theorem 2.2 
it is not the case that $\vdash_T A$. Form this it follows that $A$ is a positive part of a Hintikka formula 
by lemma 3.5. In view of Theorem 5.1, there is a model for $\mathbf{L_1}$ which falsifies the 
Hintikka formula as well as $A$. As shown in the prededing section, this model for $\mathbf{L_1}$ 
could also be the model for $\mathbf{L}$, and there $A$ is again false. But, this is contrary to the 
assumption.

This is a model-theoretic proof of Separation theorem. A syntactical proof of the original 
Separation theorem is given Takano \cite{Takano1991}.

With this we are coming to the end of the Hilbert-type axiomatic rejection for $\mathbf{L_1}$. In the 
following sections, a more constructive Gentzen-type axiomatic rejection will be developed again 
for $\mathbf{L_1}$. 

Before concluding this section, we wish to make a supplementary remark to Fundamental theorem, 
i.e. Theorem 2.1.

According to the Theorem, there obtains in a finite number of steps either a closed tableau or the 
one, which is not closed, namely open, by reducint the given formula in compliance with the stipulation as 
stated in the Theorem. Nevertheless, there might be a possibility  that some reductions are resulting in a 
closed tableau, while others do not produce any closed tableaux although starting with one and the same 
formula. We wish to show that this is not the case. In fact, if a tableau, which is open, were forthcoming 
by reducing the given formula in a way different from the successful reduction with a branch ending 
with a Hintikka formula, the Hintikka formula would be axiomatically rejected by Theorem 3.1. From this 
it follows that the $T$-transform of the formula would be not valid in first-order predicate logic with 
equality by Theorem 6.1. In view of the successful reduction, the given formula is a thesis of the 
Hilbert-type $\mathbf{L_1}$ by Theorem 2.2, and from this obtains that the $T$-transform of the 
given formula is provable in first-order predicate logic with equality, i.e. valid there by Theorem 4.1, 
namely Soundness theorem. This is a contradiction. (For this, refer to Kleene \cite{Kleene1952, 
Kleene1952a, Kleene1967}.)

Summing up the above argument, we obtain the following.
\begin{theorem}\bf 6.6 \rm (Permutability theorem) \it Once we obtain a closed tableau by reducing a formula, 
there is no possibility of getting another tableau, which is open by way of a reduction different 
from the given one. 
\end{theorem}  

\section{Gentzen-type axiomatic rejection $\mathbf{GAR}$}
The axiomatic rejection for $\mathbf{L_1}$ developed so far has been the version based upon the Hilbert-type $\mathbf{L_1}$. Thus, $\vdash_T A$, for example, has been thought of as formalized within the bounds of 
the Hilbert-type version of $\mathbf{L_1}$ notwithstanding its appearance. As mentioned earlier, this was 
also the policy adopted by \L ukasiewicz for the decision method of the Aristotelian syllogistic. 

In this and following section, we are returning to the purely Gentzen-type or tableau method version of 
$\mathbf{L_1}$ availing ourselves of its Sch\"{u}tte-style formalization as above introduced, and wish to 
develop a Gentzen-type counterpart. All the syntactical preliminaries are also understood in the 
Gentzen-style.

In the Gentzen-type axiomatic rejection, we are again starting with Theorem 2.1, namely 
Fundamental theorem. For completeness, we repeating the Theorem hereunder:

\medskip

\it 
Given a formula $($of $\mathbf{L_1})$, by reducing it by reduction rules there obtains a finite tableau, each 
branch of which ends either with a formula of th form  $F[A_+, A_-]$ or with a Hintikka formula, whereby 
a branch is extended by a reduction rule only if the formula to be reduced is not of the form  $F[A_+, A_-]$ 
and the reduction gives rise to a formula not occurring in the formula to be reduced as negative part thereof. 
\rm 

\medskip

The Theorem may be understood in the sense of the Gentzen-style formulaiton of $\mathbf{L_1}$. In 
other words, the reduction rules are the inference rules given outright in the sense of Gentzen only put
 up-side-down.

Before proceeding further, some well-known theorem will be cited of the Gentzen-type logic for 
subsequent reference: 

\begin{theorem}\bf 7.1 \rm (Thinning theorem) \it 
 
\smallskip 
$\enspace \enspace \vdash_T \enspace A \enspace \enspace \Longrightarrow \enspace \enspace 
\vdash_T \enspace F[A_+]$, 
 
$\enspace \enspace \vdash_T \enspace \sim A \enspace \Longrightarrow \enspace \enspace 
\vdash_T \enspace G[A_-]$.
\end{theorem} 
\begin{theorem}\bf 7.2 \rm (Interchange theorem) \it 

\smallskip 
$\enspace \enspace \vdash_T \enspace F[A_+, B_+] \enspace \Longrightarrow \enspace \enspace \vdash_T 
\enspace F[B_+ , A_+]$,

$\enspace \enspace \vdash_T \enspace G[A_-, B_-] \enspace \Longrightarrow \enspace \enspace \vdash_T 
\enspace G[B_- , A_-]$.
\end{theorem}
\begin{theorem}\bf 7.3 \rm (Translation theorem) \it 
 
\smallskip 
$\enspace \enspace \vdash_T \enspace F[A_+, \enspace {}_+] \enspace \Longrightarrow \enspace 
\enspace \vdash_T \enspace F[\enspace {}_+ , A_+]$,

$\enspace \enspace \vdash_T \enspace G[A_-, \enspace {}_-] \enspace \Longrightarrow \enspace 
\enspace \vdash_T \enspace G[\enspace {}_- , A_-]$,
\end{theorem}
\begin{theorem}\bf 7.4 \rm (Contraction theorem) \it 

\smallskip 
$\enspace \enspace \vdash_T \enspace F[A_+, A_+] \enspace \Longrightarrow \enspace 
\enspace \vdash_T \enspace F[A_+ , \enspace {}_+]$,

$\enspace \enspace \vdash_T \enspace G[A_-, A_-] \enspace \Longrightarrow \enspace 
\enspace \vdash_T \enspace G[A_- , \enspace {}_-]$.
\end{theorem}

All these theorems are known as structural rules \it en bloc\rm . Every one of them is easily 
proved by induction on the length of the proof (or tableau) of the assumption, except Contraction 
theorem. Hereby, it is noticed the length of the proof remains invariant or even gets shorter in 
the conclusion. 
\begin{definition}\bf 7.1 \rm  
A formula $A$ of $\mathbf{L_1}$ is \it axiomatically rejected \rm in $\mathbf{GAR}$ 
(in the Getzen-type axiomatic 
rejection) (denoted by $\dashv_T A$) if there exists a tableau of it, at least, one branch of which is 
ending with a Hintikka formula (i.e. an open tableau of it).
\end{definition}

This is a Gentzen-type counterpart of the definition of axiomatic rejection.

Roughly speaking, any Hintikka formula is now playing the role of the axioms for the Gentzen-type 
axiomatic rejection for $\mathbf{L_1}$, whereas its counterpart of the HIlbert-type $\mathbf{L_1}$, was more 
complicated as described in detail in the third section. Nevertheless, it is remembered, we come nearer the 
Gentzen-type axiomatic rejection in the same section. The theorem is to the effect ahta any formula, 
if it is rejected at all, is rejected through a Hintikka formula and after the Hintikka formula only the rule 3.13 
is made use of. (The use of Hintikka formulas as axiom for axiomatic rejection is dating from 
Inou\'{e} \cite{Inoue1989, Inoue1991} 
concerning classical propositional logic and its extensions.) 

Therefore we may regard an axiomatization for axiomatic rejection with Hintikka formulas as axioms as 
a realization of the normal form thorem in the Hilbert-type logic.

The idea of such an axiomatization for axiomatic rejection can be applied to an axiomatization for the 
set of all satisfiable formulas of classical propositional logic (and its extensions). This was pointed out 
in Inou\'{e} \cite{Inoue1993a}. 

\begin{theorem}\bf 7.5 \rm (Dichotomy theorem for the Gentzen-type $\mathbf{L_1}$) \it 
Every formula, which is not provable in the Gentzen-type $\mathbf{L_1}$ $($tableau method for 
$\mathbf{L_1})$, is axiomatically rejected in $\mathbf{GAR}$ $($in the Gentzen-type axiomatic 
rejection for $\mathbf{L_1})$.
\end{theorem}

This is a Gentzen-type counterpart of Theorem 3.2 for the Hilbert-type $\mathbf{L_1}$. The Dichotomy 
theorem is to be effect that every formula is either provable or axiomatically rejected providing us with 
a decision method for the Gentzen-type $\mathbf{L_1}$, while Theorem 3.2 itself is not yet enough to 
give a decision method for the Hilbert-type version of $\mathbf{L_1}$. This fact is a very remarkable difference between the Hilbert- and Gentzen-type versions of $\mathbf{L_1}$. 

With a view to proving Dichotomy theorem, let us assume that $A$ is not provable in the Gentzen-type 
$\mathbf{L_1}$. In view of Theorem 2.1, i.e. Fundamental theorem, there obtains by reducing $A$ a (finite) 
tableau, at least, one branch of which is ending with a Hintikka formula. $A$ is, thus, axiomatically 
rejected by Definition 7.1.
\begin{theorem}\bf 7.6 \rm (Contradiction theorem for the Gentzen-type $\mathbf{L_1}$) \it It is not the case 
that for any formula $A$ $($of  $\mathbf{L_1})$, $\vdash_T A$ and $\dashv_T A$ at the same time.
\end{theorem}
Here, $\dashv_T  A$, it is remembered, signified that $A$ is axiomatically rejected in the Gentzen-type 
$\mathbf{L_1}$. 

For the Hilbert-type version of $\mathbf{L_1}$, the Theorem was demonstrated as Theorem 6.4 in the 
preceding section. The proof there was not syntactical based upon Theorem 6.2. 

Now, we are proceeding to a syntactical proof of Theorem 7.6, i.e. Contradiction theorem.

With this in view, let us assume, if possible, $\vdash_T A$ and $\dashv_T A$ simultaneously for a 
formula $A$. 

In view of $\vdash_T A$, there obtains a closed (normal) tableau by reducing $A$. $\dashv_T A$, 
on the other hand, gives rise to a tableau, at least, one branch of which is ending with a Hintikka 
formula. This, however, is not possible on the basis of the theorem stated at the end of the sixth 
section, namely Theorem 6.4. According to the theorem, if a reduction is successful giving rise to a 
closed (normal) tableau, there is no possibility of any other reductions to fail in producing a closed 
(normal) tableau as far as they are under the proviso as stated in Theorem 2.1, namely Fundamental 
theorem. Thus, $\vdash_T A$ and $\dashv_T A$ are not compatible, and Contradiction theorem 
was proved. 

The proof of the above cited theorem given in the preceding section, i.e. Theorem 6.4, it is 
remembered, was model-theoretic resorting to Theorem 6.2. A purely syntactical proof is also 
possible of this theorem, though laborious, and we have in mind the theorem syntactically proved 
with a view to making the Gentzen-type axiomatic rejection for $\mathbf{L_1}$ purely syntactical. 

\section{Cut elimination theorem} 
In this section, we shall take up one of the highlight of this paper, namely the cut elimination theorem 
for $\mathbf{L_1}$ to be proved on the basis of the Hilbert-and Gentzen-type axiomatic rejections for 
$\mathbf{L_1}$.

As well-known, cut is a rule which is applied (in the Gentzen-Sch\"{u}tte-type formalism) is the following form (see Sch\"{u}tte \cite{Schutte1960}): 

\smallskip

8.1 $\enspace \vdash_T \enspace F[A_+], \enspace \vdash_T \enspace G[A_-] \enspace \Longrightarrow 
\enspace \enspace \vdash_T \enspace F[\enspace {}_+] \vee G[\enspace {}_-]$,

\smallskip

\noindent where $A$ is called the \it cut formula \rm of the cut-application. 
\begin{theorem}\bf 8.1 \rm (Cut elimination theorem) \it The cut rule 8.1 is a derived rule in the 
Gentzen-type version of $\mathbf{L_1}$.
\end{theorem}

In other word, cut is a rule to be dispensed with. 

Here, we shall prove cut elimination theorem model-theoretically.

With this in view, let us assume that $\vdash_T \enspace F[A_+]$ and $\vdash_T \enspace G[A_-]$. By 
Translation theorem (i.e. Theorem 7.3), they, respectively, give rise to 
$\vdash_T \enspace F[\enspace {}_+] \vee A$ and $\vdash_T \enspace G[\enspace {}_-] \vee \sim A$, 
from which we obtain, 
$$\vdash_T \enspace (F[\enspace {}_+] \vee G[\enspace {}_-] ) \vee A \enspace \mbox{and} \enspace 
\vdash_T \enspace (F[\enspace {}_+] \vee G[\enspace {}_-] ) \vee \sim A$$
by Thinning and Interchange theorems (i.e. Theorems 7.1 and 7.2). These two formulas, then, give rise to:
$$\vdash_T \enspace (F[\enspace {}_+] \vee G[\enspace {}_-] ) \vee (A \wedge \sim A)$$
$$(. \equiv . \enspace (F[\enspace {}_+] \vee G[\enspace {}_-] ) \vee 
\sim (\sim A \vee \sim \sim A)),$$
in view of its reduction by $\vee_-$. By Theorem 6.2, 
$$T(F[\enspace {}_+] \vee G[\enspace {}_-] ) \vee (A \wedge \sim A)$$
$$(\mbox{i.e.} \enspace T(F[\enspace {}_+] \vee G[\enspace {}_-] ) \vee (TA \wedge \sim TA)),$$
is valid in first-order predicate logic with equality. Since $TA \wedge \sim TA$ is contradictory, we 
have $T(F[\enspace {}_+] \vee G[\enspace {}_-] )$ is valid in the logic. Again by Theorem 6.2, we have 
the looked-for $\vdash_T \enspace F[\enspace {}_+] \vee G[\enspace {}_-]$. 

This is a semantical proof of the Cut elimination theorem for $\mathbf{L_1}$. From the standpoint of 
Gentzen-type formalization, it is not so interesting. In the sequel, we will prove it purely syntactically 
under some assumption, making use of axiomatic rejection. (In a traditional way, the cut elimination 
theorem was proved of Le\'{s}niewski' elementary) ontology by Takano \cite{Takano1991}, and 
that for $\mathbf{L_1}$ is forthcoming therefrom.) 

As well-known, with respect to other Gentzen-type logics, both the Hilbert- and Gentzen-type versions of 
$\mathbf{L_1}$ are proved to be equivalent by way of the Cut elimination theorem. The equivalence, 
it is remembered, was alsredy demonstrated in Theorem 6.2 availing ourselve of the embedding of 
$\mathbf{L_1}$ in first-order predicate logic with equality. 

Nevertheless, the equivalence thus established is confined to the provability in both versions of 
$\mathbf{L_1}$, and its counterpart for axiomatic rejection has not been proved yet. In what follows, 
we wish to demonstrate this theorem, namely:
\begin{theorem}\bf 8.2 \it For any formula $A$ of $\mathbf{L_1}$, we have 
$$\dashv_H \enspace A \enspace \Longleftrightarrow \enspace \dashv_T \enspace A.$$
\end{theorem}

The Theorem is demonstrated availing ourselves of the equivalence of the Hilbert- and Gentzen-type 
versions of $\mathbf{L_1}$.

With this in view, let us assume $\dashv_H A$, but not $\dashv_T A$. By Dichotomy theorem for the 
Gentzen-type  $\mathbf{L_1}$, we have $\vdash_T A$, which gives rise to $\vdash_H A$ by Theorem 2.2 
or the equivalence of the Hilbert-type $\mathbf{L_1}$. 

Conversely, suppose $\dashv_T A$, but not $\dashv_H A$. From this it follows in turn $\vdash_T A$ 
by Theorem 6.2. But, this is in contradiction to $\dashv_T A$ by Contradiction theorem for the 
Gentzen-type $\mathbf{L_1}$. 

This proof is evidently model-theoretic. But, we can present a syntactical proof of $\Longleftarrow$ of 
Theorem 8.2. 

That is carried out by induction on the lenght of the branch which leads to a Hintikka formula starting 
from $A$.

The basis is taken care of by Theorem 3.1 to the effect that every Hintikka formula is axiomatically 
rejected in  the Hilbert-type $\mathbf{L_1}$. 

Induction steps are dealt with by the following theses of the Hilbert-type $\mathbf{L_1}$ and induction 
hypothesis:

\smallskip

$\vdash_H \enspace G[A \vee B_-]. \supset . (G[A \vee B_-] \vee \sim A) \wedge 
(G[A \vee B_-] \vee \sim B),$

$\vdash_H \enspace G[\epsilon ab_-] \enspace \supset . \enspace G[\epsilon ab_-] 
\vee \sim \epsilon aa, $

$\vdash_H \enspace G[\epsilon ab_-, \epsilon bc_-] \enspace \supset . \enspace 
G[\epsilon ab_-, \epsilon bc_-] \vee \sim \epsilon ac, $

$\vdash_H \enspace G[\epsilon ab_-, \epsilon bb_-] \enspace \supset . \enspace 
G[\epsilon ab_-, \epsilon bb_-] \vee \sim \epsilon ba, $

$\dashv_H \enspace G[A \vee B_-] \vee \sim A \enspace $ or 
$\enspace \dashv_H \enspace G[A \vee B_-] \vee \sim B,$

$\dashv_H \enspace G[\epsilon ab_-] \vee \sim \epsilon aa, $

$\dashv_H \enspace G[\epsilon ab_-, \epsilon bc_-] \vee \sim \epsilon ac, $

$\dashv_H \enspace G[\epsilon ab_-, \epsilon bb_-] \vee \sim \epsilon ba.$

\smallskip 


Before concluding this section, we wish to present a novel syntactical proof hitherto unknown of the 
cut elimination theorem for the Genzten-type $\mathbf{L_1}$. (The proof is essentially due to the 
first author of the present paper.)

Nevertheless, an additional postulate for the Hilbert-type $\mathbf{L_1}$ is in order for demonstrating 
cut elimination theorem. It is:

\smallskip

8.4 \it No Hintikka formula of the form $A_1 \vee A_2 \vee \cdots \vee A_n$ $(n \geq 1)$ is provable 
in the Hilbert-type version of $\mathbf{L_1}$, where $A_i$ $(1 \leq i \leq n)$ is  an atomic formula or 
a negated atomic one. \rm 

\smallskip

The postulate appears to be intuitively plausible, since we can always construct a model for $\mathbf{L_1}$, 
which falsified a given Hintikka formula. From 8.4 there forthcomin th consistency of $\mathbf{L_1}$ in its 
Hilbert-type version. In fact, if the version is inconsistent, there obtains $\vdash_H \epsilon aa$, 
which contradicts that not $\vdash_H \epsilon aa$ in view of the postulate 8.4. Nevertheless, the 
consistency of $\mathbf{L_1}$ does not give rise to the postulate 8.4. This is the situation different 
from what we have in the case of classical propositional logic where the consistency is equivalent to 
the analogue of 8.4. (For the details, refer to Inou\'{e}-Ishimoto \cite{Inoue-Isimoto1992}.)

Contradiction theorem for the Hilbert-type $\mathbf{L_1}$, Theorem 6.4 it is remembered, was proved 
in the sixth section by resorting to Theorem 6.2, which was model-theoretic. 

A purely syntactic proof of the Contradiction theorem for the Hilbert-type version of 
$\mathbf{L_1}$ is carried out by induction on the length of the axiomatic rejection of the given $A$, 
where the postulate 8.4 is playing an important role.

As easily seen, for proving Contradiction theorem it is sufficient to derive not $\vdash_H A$ from 
$\dashv_H A$. 

The basis cases (3.11 and 3.12) are taken care of by 8.4. Indeed, $\vdash_H \epsilon aa$ and 
$\dashv_H \epsilon aa$ ($\vdash_H \sim \epsilon aa$ and $\dashv_H \sim \epsilon aa$) are not the 
case simultaneously in view of the postulate 8.4.

For taking care of induction steps, suppose $\dashv_H \epsilon aa$ is obtained from 
$\vdash_H A \supset B$ and $\dashv_H B$ by the rule 3.13. Suppose, further, $\dashv_H A$, which 
gives rise to $\vdash_H B$ by detachment against induction hypothesis. This takes care of the rule 3.13 for axiomatic rejection. 

With a view to dealing with the rule 3.14, suppose $\dashv_H B$ and $B$ is forthcoming from $A$ by a 
uniform substitution for some name variables occurring in $A$. Further, suppose, if possible, that 
$\vdash_H B$, which gives rise to  $\vdash_H A$ against induction hypothesis. 

Lastly, assume $\dashv_H A \vee \epsilon aa$ is obtained from $\dashv_H A$ by means of the 
rule 3.15,  $A \vee \epsilon aa$ is a Hintikka formula with the condition in 3.15,  
$\dashv_H A \vee \epsilon aa$ does not hold by the postulate 8.4. 

This completes the syntactical proof of the Contradiction theorem for the Hilbert-type  version of 
$\mathbf{L_1}$ (under the postulate 8.4). 

We are now in a position to demonstrate the Cut elimination theorem for the Genzten-type 
$\mathbf{L_1}$ on the basis of the Contradiction for the Hilbert-type $\mathbf{L_1}$ just proved.

The demonstration of the Cut elimination theorem prodeeds in the following way: 

To start with, suppose $\vdash_T F[A_+]$ and $\vdash_T G[A_-]$, which, respectively, give rise to 
$\vdash_H F[A_+]$ and $\vdash_H G[A_-]$ by Theorem 2.2. They will be referred to as (*) and (**) 
below, respectively. 

With a view to obtaining $\vdash_T \enspace F[\enspace {}_+] \vee G[\enspace {}_-]$, let us assume, 
if possible, to the contrary, namely not $\vdash_T \enspace F[\enspace {}_+] \vee G[\enspace {}_-]$, 
which in turn yields $\dashv_T \enspace F[\enspace {}_+] \vee G[\enspace {}_-]$ by the Dichotomy 
theorem for Gentzen-type $\mathbf{L_1}$, namely Theorem 7.5. In view of Definition 7.1, there 
obtains a Hintikka formula by reducing $F[\enspace {}_+] \vee G[\enspace {}_-]$. By Theorem 3.1 and 
Lemma 3.5, we, then, have $\dashv_H F[\enspace {}_+] \vee G[\enspace {}_-]$. 

We, now, wish to derived $\vdash_H F[\enspace {}_+] \vee G[\enspace {}_-]$ on the basis of (*) and 
(**), which is going on in the following way: 

\smallskip

(1)  $\enspace \vdash_H \enspace F[A_+]$ $\enspace \enspace \enspace \enspace \enspace$ (*), 

(2)  $\enspace \vdash_H \enspace G[A_-]$ $\enspace \enspace \enspace \enspace \enspace$ (**), 

(3)  $\enspace \vdash_H \enspace F[\enspace {}_+] \vee A$ $\enspace \enspace \enspace 
\enspace \enspace$ (1), Lemma 2.1, 

(4)  $\enspace \vdash_H \enspace G[\enspace {}_-] \vee A$ $\enspace \enspace \enspace 
\enspace \enspace$ (2), Lemma 2.1, 

(5)  $\enspace \vdash_H \enspace (F[\enspace {}_+] \vee G[\enspace {}_-]) \vee A$ $\enspace \enspace \enspace \enspace \enspace$ (3), tautology, 

(6)  $\enspace \vdash_H \enspace  (F[\enspace {}_+] \vee G[\enspace {}_-]) \vee A$ $\enspace \enspace \enspace 
\enspace \enspace$ (4), tautology, 

(7)  $\enspace \vdash_H \enspace (F[\enspace {}_+] \vee G[\enspace {}_-]) \vee (A \wedge \sim A)$ 
$\enspace \enspace \enspace \enspace \enspace$ (5), (6), tautology, 

(8)  $\enspace \vdash_H \enspace F[\enspace {}_+] \vee G[\enspace {}_-] . \supset . 
F[\enspace {}_+] \vee G[\enspace {}_-]$ $\enspace \enspace \enspace 
\enspace \enspace$ tautology, 

(9)  $\enspace \vdash_H \enspace A \wedge \sim A . \supset . 
F[\enspace {}_+] \vee G[\enspace {}_-]$ $\enspace \enspace \enspace 
\enspace \enspace$ tautology, 

(10)  $\enspace \vdash_H \enspace 
F[\enspace {}_+] \vee G[\enspace {}_-]$ $\enspace \enspace \enspace 
\enspace \enspace$ (7), (8), (9), tautology.

\smallskip

The last formula (10), namely $\vdash_H \enspace F[\enspace {}_+] \vee G[\enspace {}_-]$, 
however contradicts $\dashv_H \enspace F[\enspace {}_+] \vee G[\enspace {}_-]$ as above 
obtained in view of the Contradiction theorem for the Hilbert-type $\mathbf{L_1}$.

This completes the demonstration of Cut elimination theorem on the basis of the postulate 8.4.

Here, it is remarked in passing that Theorem 6.2, in virture of which we are allowed to obtain 
$\vdash_T A$ from $\vdash_H A$ is not employing Cut elimination theorem for its proof. 
Otherwise, it would be preposterous. 

This kind of proof of Cut elimination theorem was first explored in Inou\'{e}-Ishimoto 
\cite{Inoue-Isimoto1992} for 
classical propositional logic and will be made use of in other syllogistic systems as shown in 
Inou\'{e}-Ishimoto \cite{Inoue-Ishimoto-fc}.

By the Cut elimination theorem for $\mathbf{L_1}$, we obtain, in a routine way, a syntactical proof of: 

\smallskip 

8.5 $\enspace \enspace \vdash_H \enspace A \enspace \Longrightarrow 
\enspace \enspace \vdash_T \enspace A$, 

\smallskip 

\noindent which was model-theoretically proved in Theorem 6.2. 

As a concluding remark of the present section, we shall show the following.
\begin{theorem} \bf 8.3 \it 
The following statements are equivalent:

\smallskip

$(1)$ The Cut elimination theorem for $\mathbf{L_1}$ holds,

$(2)$  No Hintikka formula of the form $A_1 \vee A_2 \vee \cdots \vee A_n$ $(n \geq 1)$ is provable 
in the Hilbert-type version of $\mathbf{L_1}$, where $A_i$ $(1 \leq i \leq n)$ is  an atomic formula or 
a negated atomic one. $($= the postulate 8.4$)$, 

$(3)$ $\enspace \enspace \vdash_H \enspace A \enspace \Longrightarrow 
\enspace \enspace \vdash_T \enspace A$, 

$(4)$) Contradiction and Dichotomy theorems for the Hilbert-type $\mathbf{L_1}$ hold.
\end{theorem}

We have already demonstrated the proof of (2) $\Rightarrow$ (4) $\Rightarrow$ (1) above. And we 
mentioned (1) $\Rightarrow$ (3) in 8.5. We shall prove the implication (3) $\Rightarrow$ (2). 
The proof of it is not difficult to think of. Suppose 
$$A_1\vee A_2 \vee \dots \vee A_n \enspace \enspace (n \geq 1)$$
is a Hintikka formula, where  $A_i$ ($1 \leq i \leq n$) is an atomic formula or a negated atomic one. Thus, 
not $\vdash_TA_1\vee A_2 \vee \dots \vee A_n$. By the contraposition of (3), the formula is not 
provable in the Hilbert-type $\mathbf{L_1}$. 

This completes the proof of Theorem 8.3. We wish to emphasize that the proof is purely syntactical.

\section{Characterization theorem and  axiomatic rejection with Hintikka formulas as axioms}

The reader will find a similar equivalence as Theorem 8.3 in Inou\'{e}-Ishimoto \cite{Inoue-Isimoto1992} for classical 
propositional logic, on the basis of which the argument for Theorem 8.3 was developed. Such 
equivalences would hold for a variety of logics, if (2) is appropriately changed for a given logic. 


Inou\'{e} \cite{Inoue1995a} obtained the following theorem.

\begin{theorem}\bf 9.1 \rm (Theorem 1.1 in \cite{Inoue1995a}) \it 
The following statements are equivalent: 

\smallskip

$(1)$ No Hintikka formula is provable in the Hilbert-type $\mathbf{L_1}$\footnote{Please let the first author of this paper give a correction of his paper. The original statement of Theorem 1.1 (i) in \cite{Inoue1995a} is 
not correct. `$\mathbf{HL_1}$' of Theorem 1.1 (i) in \cite{Inoue1995a} should be `the Hilbert-type 
$\mathbf{L_1}$.}, 

$(2)$ The Hilbert-type $\mathbf{L_1}$ is \L -decidable with respect to $\mathbf{HL_1}$ $($i.e. the set of all formula of 
$\mathbf{L_1}$ is the disjoint union of the set of all the theorem of the Hilbert-type $\mathbf{L_1}$ and that of 
all the theorem of $\mathbf{HL_1})$,

$(3)$ The Cut elimination theorem for the Gentzen-type $\mathbf{L_1}$ holds,

$(4)$ For any formula $A$ of $\mathbf{L_1}$, if $A$ is provable in the Hilbert-type $\mathbf{L_1}$, then 
it is provable in the the Gentzen-type $\mathbf{L_1}$.

\smallskip

\end{theorem}

Theorem 9.1 was proved on the basis of the Hilbert-type axiomatic rejection $\mathbf{HL_1}$ for 
$\mathbf{L_1}$ which is defined with Hintikka formulas as axioms in \cite{Inoue1995a} in the 
following.

Axioms:

\noindent 9.1 For every Hintikka formula $A$ of $\mathbf{L_1}$, $\enspace \vdash_{\mathbf{HL_1}} 
\enspace A$, 

Rule:

\noindent 9.2 $\enspace \enspace \vdash_H \enspace A \supset B, \enspace 
\vdash_{\mathbf{HL_1}} \enspace B 
\Longrightarrow \enspace \enspace \vdash_{\mathbf{HL_1}} \enspace A$,
\medskip

\noindent where $\vdash_{\mathbf{HL_1}} A$ means that $A$ is axiomatically rejected by 
$\mathbf{HL_1}$.

Combining with Theorem 8.3 with the above Theorem 9.1, we obtain our last principal result of this paper 
as follows.

\begin{theorem}\bf 9.2 \rm (Characterization Theorem) \it 
The following statements are equivalent:

\smallskip

$(1)$ The Cut elimination theorem for the Gentzen-type $\mathbf{L_1}$ $($tableau method$)$ holds,

$(2)$  No Hintikka formula of the form $A_1 \vee A_2 \vee \cdots \vee A_n$ $(n \geq 1)$ is provable 
in the Hilbert-type version of $\mathbf{L_1}$, where $A_i$ $(1 \leq i \leq n)$ is an atomic formula or 
a negated atomic one,

$(3)$ For any formula $A$ of $\mathbf{L_1}$, if $A$ is provable in the Hilbert-type $\mathbf{L_1}$, then 
it is provable in the the Gentzen-type $\mathbf{L_1}$ $($tableau method$)$, 

$(4)$ Contradiction and Dichotomy theorems for the Hilbert-type $\mathbf{L_1}$ hold.

$(5)$ No Hintikka formula is provable in the Hilbert-type $\mathbf{L_1}$,

$(6)$ The Hilbert-type $\mathbf{L_1}$ is \L -decidable with respect to $\mathbf{HL_1}$ $($i.e. the set of all formula 
of the Hilbert-type $\mathbf{L_1}$ is the disjoint union of the set of all the theorem of $\mathbf{L_1}$ 
and that of all the theorem of $\mathbf{HL_1)}$. 

\smallskip

\end{theorem}

From Theorem 9.2.(4), (6), Theorems 3.2 and 6.4, we have 

\begin{theorem}\bf 9.3 \it For any formula $A$ of $\mathbf{L_1}$, 
$$\vdash_{\mathbf{HL_1}} A \enspace \Longleftrightarrow \enspace \vdash_{\mathbf{HAR}} A.$$
\end{theorem}

Thus this means that the $\mathbf{HL_1}$ has the same strength with $\mathbf{HAR}$.

%

%
%
%

$$ $$

\centerline{\large \bf Appendix by Arata Ishimoto, the second author \rm}

\medskip 

\noindent Up to the last section\footnote{In the original version of this paper, 
this appendix is the section 9.}, we have mainly been concerned with the technical matters of $\mathbf{L_1}$, 
namely the propositional fragment of Le\'{s}niewski's ontology ignoring its philosophical implications. 
But, Le\'{s}niewski's ontology is a system intended to be a logic in the wider sense of the word, 
not a mere formalism as emphasized by Le\'{s}niewski himself. Roughly speaking, Le\'{s}niewski's 
ontology is an ontology in the traditional sense of the word. This is to the effect that Le\'{s}niewski's 
ontology is a science to inquire into the most general aspects of the entities existing in the 
world.

From such a poinf of view, we wish to scrutinize some philosphico-ontological problem underlying 
the logical techniques developed so far. More specifically, we are looking into the 
philosophico-ontological significance of $\mathbf{L_1}$ to be seen under formalism. This, it is 
believed, is the very task of philosophical logic.

As well-known, Le\'{s}niewski's ontology has traditionally been interpreted in the spirit a rather 
extreme nominalism, which has culminated in the so-called `reism' as propounded by 
Kotarbi\'{n}ski. Reism is a philosophy which advocates an ontology that only material things are 
legitimate intities in existence, and beyond them there is nothing. (For reism, refer to 
Wole\'{n}ski \cite{Wolenski1990} and also to Sinisi \cite[p. 59]{Sinisi1983} with respect to 
Le\'{s}niewski's ontology.) 

Nevertheless, we wish to oppose to such a nominalism another ontology diagonally different 
therefrom. To be more specific, a conceptula realism will be advanced here as an alternative 
ontology to underlie Le\'{s}niewski's ontology or its propositional fragment $\mathbf{L_1}$. The 
conceptual realism we are defending is a Platonist philosophy to the effect that only concepts or 
universals are in existence independently and they are remaining invariant through every interpretation. 
In reference to the $\mathbf{L_1}$ we have developed in this paper, name variables are 
representing concepts, i.e. universals, and they are susceptible of a large number of different 
interpretations. This was shown in the concrete in the fifth and sixth sections with respect to the 
model construction for $\mathbf{L_1}$. More specifically, Theorem 6.2 tells us that $A$ is a thesis of 
$\mathbf{L_1}$ if and only if $TA$ is valid in first-order predicate logic with equality, where $T$ is 
a translation to transform a formula of $\mathbf{L_1}$ into its correspondent in predicate logic. This 
is an embedding theorem of $\mathbf{L_1}$ in first-order predicate logic. As shown in the fifth section, 
there are a variety of possibilities of defining moels for $\mathbf{L_1}$. Thus, for one and the same 
$A$, we have a large number of different interpretations, and every one of them  makes $TA$ true. 
As indicated in the sixth section, $A$ is also made true in the models which are defined in terms of 
those for $TA$.

In defining a model for the given $A$, it is remarked, each name variable involved in $A$ is assigned a set of individuals (in the sense of predicate logic), and the set varies from one interpretation to another. 
Nevertheless, the name variable, on the basis of which we are defining sets of individuals, is trascending 
all these sets remaining the same concept or universal.

This is nothing but the conceptual realism we are defending availing ourselves of the technical 
aparatus as developed up to the preceding section. Put it the other way round, $\mathbf{L_1}$ is 
a logic deprived of individuals, namely the entities belonging to the lowest type which are called 
out temporarily legitimate entities in the proposed ontology of conceptual realism. (For such a 
conceptual realism, refer also to Ishimoto \cite{Ishimoto1986}.)





\bigskip 

\bigskip 

\noindent \bf Acknowledgments. \rm The first author of this paper, as the representative of us, would like to 
thank the late Professor V. A. Smirnov for inviting us to Institute of Philosophy, Russian Academy of Sciences in 
Moscow in order to present this work at the conference. This paper is an enlarged and refined version of the 
paper of it. For Professor Smirnov, see Karpenko \cite{Karpenko2000}, Bystrov \cite{Bystrov1996} and 
Finn \cite{Finn2000}. The first author would also like to thank the late Professor Emeritus Arata 
Ishimoto, my teacher, and Mr. Mitsunori Kobayashi, my research friend, for fruitful research and discussions.



\noindent Takao Inou\'{e}

\noindent 1. Meiji Pharmaceutical University \\ Department of Medical Molecular Informatics \\Tokyo, Japan 
\bigskip

\noindent 2. Hosei University

\noindent Graduate School of Science and Engineering\\ Tokyo, Japan
\bigskip

\noindent 3. Hosei University\\ Faculty of Science and Engineering\\ Department of Applied Informatics\\ Tokyo, Japan
\medskip

\noindent ta-inoue@my-pharm.ac.jp \\ takao.inoue.22@hosei.ac.jp \\ takaoapple@gmail.com
$$ $$
\noindent Arata Ishimoto

\noindent Professor Emeritus of Tokyo Institute of Technology \\ Tokyo, Japan \\ Deceased 
$$ $$
\noindent Mitsunori Kobayashi

\noindent Logician and Composer, Japan

\end{document}